\documentclass{article}
\usepackage{amsmath,amsfonts,amssymb,bm,hyperref}
\numberwithin{equation}{section}
\newtheorem{theorem}{Theorem}[section]
\newtheorem{lemma}{Lemma}[section]

\newtheorem{remark}{Remark}[section]
\newtheorem{definition}{Definition}[section]
\usepackage[T1]{fontenc}
\usepackage[utf8]{inputenc}
\usepackage{authblk}

\author{Olga Chervova\thanks{O.Chervova@ucl.ac.uk, \url{http://www.homepages.ucl.ac.uk/\~ucahoch/}} }
\author{Robert J.~Downes\thanks{R.Downes@ucl.ac.uk, \url{http://www.homepages.ucl.ac.uk/\~zcahc37/}} }
\author{Dmitri Vassiliev\thanks{D.Vassiliev@ucl.ac.uk, \url{http://www.homepages.ucl.ac.uk/\~ucahdva/}}}
\affil{Department of Mathematics,
University College London,\\ Gower Street, London WC1E 6BT, UK}

\begin{document}

\title{The spectral function of a first order elliptic system}
\maketitle
\begin{abstract}
We consider an elliptic self-adjoint first order pseudodifferential
operator acting on columns of complex-valued half-densities over
a connected compact manifold without boundary. The eigenvalues of the principal
symbol are assumed to be simple but no assumptions are made on their
sign, so the operator is not necessarily semi-bounded.
We study the following objects:
the propagator (time-dependent operator which solves
the Cauchy problem for the dynamic equation),
the spectral function (sum of squares of Euclidean norms
of eigenfunctions evaluated at a given point of the manifold,
with summation carried out over all eigenvalues between zero
and a positive~$\lambda$) and
the counting function (number of eigenvalues between zero
and a positive~$\lambda$).
We derive explicit two-term asymptotic formulae for all three.
For the propagator ``asymptotic'' is understood as asymptotic
in terms of smoothness, whereas for the spectral and counting
functions ``asymptotic'' is understood as asymptotic with respect
to $\lambda\to+\infty$.
\end{abstract}


\



\section{Main results}
\label{Main results}

The aim of the paper is to extend the classical results of \cite{DuiGui}
to systems. We are motivated by the observation that,
to our knowledge, all previous publications on systems give
formulae for the second asymptotic coefficient that are either
incorrect or incomplete (i.e.~an algorithm for the calculation
of the second asymptotic coefficient rather than an actual formula).
The appropriate bibliographic review is presented in
Section~\ref{Bibliographic review}.

Consider a first order classical pseudodifferential
operator $A$ acting on columns
$v=\begin{pmatrix}v_1&\ldots&v_m\end{pmatrix}^T$
of complex-valued half-densities
over a connected compact $n$-dimensional manifold $M$ without boundary.
Throughout this paper we assume that $m,n\ge2\,$.

We assume the coefficients of the operator $A$ to be infinitely smooth. We also
assume that the operator $A$ is formally self-adjoint (symmetric):
$\int_Mw^*Av\,dx=\int_M(Aw)^*v\,dx$ for all infinitely smooth
$v,w:M\to\mathbb{C}^m$. Here and further on
the superscript $\,{}^*\,$ in matrices, rows and columns
indicates Hermitian conjugation in $\mathbb{C}^m$
and $dx:=dx^1\ldots dx^n$, where $x=(x^1,\ldots,x^n)$ are local
coordinates on $M$.

Let $A_1(x,\xi)$ be the principal symbol of the operator $A$.
Here $\xi=(\xi_1,\ldots,\xi_n)$ is the variable dual to the position
variable $x$; in physics literature the $\xi$ would be referred to
as \emph{momentum}. Our principal symbol $A_1$ is an $m\times m$
Hermitian matrix-function on $T'M:=T^*M\setminus\{\xi=0\}$,
i.e.~on the cotangent bundle with the zero section removed.

Let $h^{(j)}(x,\xi)$ be the eigenvalues of the principal symbol. We
assume these eigenvalues to be nonzero (this is a version of the
ellipticity condition) but do not make any assumptions on their
sign. We also assume that the eigenvalues $h^{(j)}(x,\xi)$ are
simple for all $(x,\xi)\in T'M$. The techniques developed in our
paper do not work in the case when eigenvalues of the principal
symbol have variable multiplicity, though they could probably be
adapted to the case of constant multiplicity different from
multiplicity 1. The use of the letter ``$h$'' for an eigenvalue of
the principal symbol is motivated by the fact that later it will
take on the role of a Hamiltonian, see formula (\ref{Hamiltonian system of equations}).

We enumerate the eigenvalues of the principal symbol
$h^{(j)}(x,\xi)$ in increasing order, using a positive index
$j=1,\ldots,m^+$ for positive $h^{(j)}(x,\xi)$ and a negative index
$j=-1,\ldots,-m^-$ for negative $h^{(j)}(x,\xi)$. Here $m^+$ is the
number of positive eigenvalues of the principal symbol and $m^-$ is
the number of negative ones. Of course, $m^++m^-=m$.

Under the above assumptions $A$ is a self-adjoint operator, in the
full functional analytic sense, in the Hilbert space
$L^2(M;\mathbb{C}^m)$ (Hilbert space of square integrable
complex-valued column ``functions'') with domain $H^1(M;\mathbb{C}^m)$
(Sobolev space of complex-valued column ``functions'' which are
square integrable together with their first partial derivatives) and
the spectrum of $A$ is discrete. These facts are easily established
by constructing the parametrix (approximate inverse) of the operator
$A+iI$.

Let $\lambda_k$ and
$v_k=\begin{pmatrix}v_{k1}(x)&\ldots&v_{km}(x)\end{pmatrix}^T$ be
the eigenvalues and eigenfunctions of the operator $A$. The
eigenvalues $\lambda_k$ are enumerated in increasing order with
account of multiplicity,
using a positive index $k=1,2,\ldots$ for positive $\lambda_k$
and a nonpositive index $k=0,-1,-2,\ldots$ for nonpositive $\lambda_k$.
If the operator $A$ is bounded from below (i.e.~if $m^-=0$)
then the index $k$ runs from some integer value to $+\infty$;
if the operator $A$ is bounded from above (i.e.~if $m^+=0$)
then the index $k$ runs from $-\infty$ to some integer value;
and if the operator $A$ is unbounded from above
and from below (i.e.~if $m^+\ne0$ and $m^-\ne0$)
then the index $k$ runs from $-\infty$ to $+\infty$.

\

We will be studying the following three objects.

\

\textbf{Object 1.}
Our first object of study is the \emph{propagator},
which is the one-parameter family of operators defined as
\begin{equation}
\label{definition of wave group}
U(t):=e^{-itA}
=\sum_k e^{-it\lambda_k}v_k(x)\int_M[v_k(y)]^*(\,\cdot\,)\,dy\,,
\end{equation}
$t\in\mathbb{R}$.
The propagator provides a solution to the Cauchy problem
\begin{equation}
\label{initial condition most basic}
\left.w\right|_{t=0}=v
\end{equation}
for the dynamic equation
\begin{equation}
\label{dynamic equation most basic}
D_tw+Aw=0\,,
\end{equation}
where $D_t:=-i\partial/\partial t$.
Namely, it is easy to see that if the column of half-densities $v=v(x)$
is infinitely smooth,
then, setting
$\,w:=U(t)\,v$, we get a time-dependent column of half-densities $w(t,x)$
which is also infinitely smooth
and which satisfies the equation
(\ref{dynamic equation most basic})
and the initial condition
(\ref{initial condition most basic}).
The use of the letter ``$U$'' for the propagator is motivated by the
fact that for each $t$ the operator $U(t)$ is unitary.

\

\textbf{Object 2.}
Our second object of study is the \emph{spectral function},
which is the real density defined as
\begin{equation}
\label{definition of spectral function}
e(\lambda,x,x):=\sum_{0<\lambda_k<\lambda}\|v_k(x)\|^2,
\end{equation}
where $\|v_k(x)\|^2:=[v_k(x)]^*v_k(x)$ is the square of the
Euclidean norm of the eigenfunction $v_k$ evaluated at the point
$x\in M$ and $\lambda$ is a positive parameter (spectral parameter).

\

\textbf{Object 3.}
Our third and final object of study is the \emph{counting function}
\begin{equation}
\label{definition of counting function}
N(\lambda):=\,\sum_{0<\lambda_k<\lambda}1\ =\int_Me(\lambda,x,x)\,dx\,.
\end{equation}
In other words, $N(\lambda)$ is the number of eigenvalues $\lambda_k$
between zero and $\lambda$.

\

It is natural to ask the question: why, in defining the spectral function
(\ref{definition of spectral function})
and the counting function
(\ref{definition of counting function}),
did we choose to perform summation
over all \emph{positive} eigenvalues up to a given positive $\lambda$
rather than
over all \emph{negative} eigenvalues up to a given negative $\lambda$?
There is no particular reason. One case reduces to the other by the change
of operator $A\mapsto-A$. This issue will be revisited in
Section~\ref{Spectral asymmetry}.

Further on we assume that $m^+>0$, i.e.~that the operator
$A$ is unbounded from above.

\

Our objectives are as follows.

\

\textbf{Objective 1.}
We aim to construct the propagator
(\ref{definition of wave group})
explicitly in terms of
oscillatory integrals, modulo an integral operator with an
infinitely smooth, in the variables $t$, $x$ and $y$, integral kernel.

\

\textbf{Objectives 2 and 3.}
We aim to derive, under appropriate assumptions on Hamiltonian
trajectories, two-term asymptotics for the spectral function
(\ref{definition of spectral function})
and the counting function
(\ref{definition of counting function}),
i.e.~formulae of the type
\begin{equation}
\label{two-term asymptotic formula for spectral function}
e(\lambda,x,x)=a(x)\,\lambda^n+b(x)\,\lambda^{n-1}+o(\lambda^{n-1}),
\end{equation}
\begin{equation}
\label{two-term asymptotic formula for counting function}
N(\lambda)=a\lambda^n+b\lambda^{n-1}+o(\lambda^{n-1})
\end{equation}
as $\lambda\to+\infty$.
Obviously, here we expect the real constants $a$, $b$ and real densities
$a(x)$, $b(x)$ to be related in accordance with
\begin{equation}
\label{a via a(x)}
a=\int_Ma(x)\,dx,
\end{equation}
\begin{equation}
\label{b via b(x)}
b=\int_Mb(x)\,dx.
\end{equation}

\

It is well known that the above three objectives are closely
related: if one achieves Objective 1, then Objectives 2 and 3 follow via
Fourier Tauberian theorems \cite{DuiGui,mybook,ivrii_book,Safarov_Tauberian_Theorems}.

\

We are now in a position to state our main results.

\

\textbf{Result 1.}
We construct the propagator as a sum of $m$ oscillatory integrals
\begin{equation}
\label{wave group as a sum of oscillatory integrals}
U(t)\overset{\operatorname{mod}C^\infty}=\sum_j
U^{(j)}(t)\,,
\end{equation}
where the phase function of each oscillatory integral
$U^{(j)}(t)$ is associated with the corresponding
Hamiltonian $h^{(j)}(x,\xi)$. The symbol of the oscillatory integral
$U^{(j)}(t)$ is a complex-valued $m\times m$ matrix-function
$u^{(j)}(t;y,\eta)$, where $y=(y^1,\ldots,y^n)$ is the position of the
source of the wave (i.e.~this is the same $y$ that appears in formula
(\ref{definition of wave group})) and
$\eta=(\eta_1,\ldots,\eta_n)$ is the corresponding dual variable
(covector at the point $y$).
When $|\eta|\to+\infty$, the symbol admits an asymptotic expansion
\begin{equation}
\label{decomposition of symbol of OI into homogeneous components}
u^{(j)}(t;y,\eta)=u^{(j)}_0(t;y,\eta)+u^{(j)}_{-1}(t;y,\eta)+\ldots
\end{equation}
into components positively homogeneous in $\eta$, with the subscript
indicating degree of homogeneity.

The formula for the principal symbol of the oscillatory integral
$U^{(j)}(t)$ is known
\cite{SafarovDSc,NicollPhD}
and reads as follows:
\begin{multline}
\label{formula for principal symbol of oscillatory integral}
u^{(j)}_0(t;y,\eta)=
[v^{(j)}(x^{(j)}(t;y,\eta),\xi^{(j)}(t;y,\eta))]
\,[v^{(j)}(y,\eta)]^*
\\
\times\exp
\left(
-i\int_0^tq^{(j)}(x^{(j)}(\tau;y,\eta),\xi^{(j)}(\tau;y,\eta))\,d\tau
\right),
\end{multline}
where $v^{(j)}(z,\zeta)$ is the normalised eigenvector of the principal
symbol $A_1(z,\zeta)$ corresponding to the eigenvalue (Hamiltonian)
$h^{(j)}(z,\zeta)$,
\ $(x^{(j)}(t;y,\eta),\xi^{(j)}(t;y,\eta))$ is the Hamiltonian trajectory
originating from the point $(y,\eta)$, i.e.~solution of the system of
ordinary differential equations (the dot denotes differentiation in $t$)
\begin{equation}
\label{Hamiltonian system of equations}
\dot x^{(j)}=h^{(j)}_\xi(x^{(j)},\xi^{(j)}),
\qquad
\dot\xi^{(j)}=-h^{(j)}_x(x^{(j)},\xi^{(j)})
\end{equation}
subject to the initial condition $\left.(x^{(j)},\xi^{(j)})\right|_{t=0}=(y,\eta)$,
\ $q^{(j)}:T'M\to\mathbb{R}$ is the function
\begin{equation}
\label{phase appearing in principal symbol}
q^{(j)}:=[v^{(j)}]^*A_\mathrm{sub}v^{(j)}
-\frac i2
\{
[v^{(j)}]^*,A_1-h^{(j)},v^{(j)}
\}
-i[v^{(j)}]^*\{v^{(j)},h^{(j)}\}
\end{equation}
and
\begin{equation}
\label{definition of subprincipal symbol}
A_\mathrm{sub}(z,\zeta):=
A_0(z,\zeta)+\frac i2
(A_1)_{z^\alpha\zeta_\alpha}(z,\zeta)
\end{equation}
is the subprincipal symbol of the operator $A$,
with the subscripts $z^\alpha$ and $\zeta_\alpha$
indicating partial derivatives and
the repeated index $\alpha$ indicating summation over $\alpha=1,\ldots,n$.
Curly brackets in formula
(\ref{phase appearing in principal symbol})
denote the Poisson bracket on matrix-functions
\begin{equation}
\label{Poisson bracket on matrix-functions}
\{P,R\}:=P_{z^\alpha}R_{\zeta_\alpha}-P_{\zeta_\alpha}R_{z^\alpha}
\end{equation}
and its further generalisation
\begin{equation}
\label{generalised Poisson bracket on matrix-functions}
\{P,Q,R\}:=P_{z^\alpha}QR_{\zeta_\alpha}-P_{\zeta_\alpha}QR_{z^\alpha}\,.
\end{equation}

As the derivation of formula
(\ref{formula for principal symbol of oscillatory integral})
was previously performed only in theses \cite{SafarovDSc,NicollPhD},
we repeat it in Sections
\ref{Algorithm for the construction of the wave group}
and
\ref{Leading transport equations}
of our paper.
Our derivation differs slightly from that in \cite{SafarovDSc} and \cite{NicollPhD}.

Formula (\ref{formula for principal symbol of oscillatory integral})
is invariant under changes of local coordinates on the manifold $M$,
i.e.~elements of the $m\times m$ matrix-function
$u^{(j)}_0(t;y,\eta)$ are scalars on
$\mathbb{R}\times T'M$.
Moreover, formula (\ref{formula for principal symbol of oscillatory integral})
is invariant under the transformation of the eigenvector of
the principal symbol
\begin{equation}
\label{gauge transformation of the eigenvector}
v^{(j)}\mapsto e^{i\phi^{(j)}}v^{(j)},
\end{equation}
where
\begin{equation}
\label{phase appearing in gauge transformation}
\phi^{(j)}:T'M\to\mathbb{R}
\end{equation}
is an arbitrary smooth function.
When some quantity is defined up to the action of a certain
transformation, theoretical physicists refer to such a
transformation as a \emph{gauge transformation}. We follow this
tradition.
Note that our particular gauge
transformation
(\ref{gauge transformation of the eigenvector}),
(\ref{phase appearing in gauge transformation})
is quite common in quantum mechanics:
when $\phi^{(j)}$ is a function of the position variable $x$ only
(i.e.~when $\phi^{(j)}:M\to\mathbb{R}$) this gauge transformation
is associated with electromagnetism.

Both Y.~Safarov \cite{SafarovDSc} and W.J.~Nicoll \cite{NicollPhD} assumed
that the operator $A$ is semi-bounded from below but this assumption
is not essential and their formula
(\ref{formula for principal symbol of oscillatory integral})
remains true in the more general case that we are dealing with.

However, knowing the principal symbol
(\ref{formula for principal symbol of oscillatory integral})
of the oscillatory integral $U^{(j)}(t)$
is not enough if one wants to derive two-term
asymptotics
(\ref{two-term asymptotic formula for spectral function})
and
(\ref{two-term asymptotic formula for counting function}).
One needs information about $u^{(j)}_{-1}(t;y,\eta)$,
the component of the symbol of the
oscillatory integral $U^{(j)}(t)$
which is positively homogeneous in $\eta$ of degree~-1,
see formula (\ref{decomposition of symbol of OI into homogeneous components}),
but here the problem is that $u^{(j)}_{-1}(t;y,\eta)$
is not a true invariant in the sense that it depends on the choice of
phase function in the oscillatory integral. We overcome this difficulty
by observing that $U^{(j)}(0)$ is a pseudodifferential operator, hence,
it has a well-defined subprincipal symbol
$[U^{(j)}(0)]_\mathrm{sub}$. We prove that
\begin{equation}
\label{subprincipal symbol of OI at time zero}
\operatorname{tr}[U^{(j)}(0)]_\mathrm{sub}
=-i\{[v^{(j)}]^*,v^{(j)}\}
\end{equation}
and subsequently show that information contained in formulae
(\ref{formula for principal symbol of oscillatory integral})
and
(\ref{subprincipal symbol of OI at time zero})
is sufficient for the derivation of two-term
asymptotics
(\ref{two-term asymptotic formula for spectral function})
and
(\ref{two-term asymptotic formula for counting function}).

Note that the RHS of formula (\ref{subprincipal symbol of OI at time zero})
is invariant under the gauge transformation
(\ref{gauge transformation of the eigenvector}),
(\ref{phase appearing in gauge transformation}).

Formula~(\ref{subprincipal symbol of OI at time zero})
plays a central role in our paper.
Sections~\ref{Algorithm for the construction of the wave group}
and~\ref{Leading transport equations}
provide auxiliary material needed for the proof of
formula~(\ref{subprincipal symbol of OI at time zero}),
whereas the actual proof
of formula~(\ref{subprincipal symbol of OI at time zero})
is given in Section~\ref{Proof of formula}.

Let us elaborate briefly on the geometric meaning of the
RHS of (\ref{subprincipal symbol of OI at time zero})
(a more detailed exposition is presented in Section~\ref{U(1) connection}).
The eigenvector of the principal
symbol is defined up to a gauge transformation
(\ref{gauge transformation of the eigenvector}),
(\ref{phase appearing in gauge transformation})
so it is natural to introduce a $\mathrm{U}(1)$ connection on $T'M$ as
follows: when parallel transporting an eigenvector of the principal
symbol along a curve in $T'M$ we require that the derivative of the
eigenvector along the curve be orthogonal to the eigenvector itself.
This is equivalent to the introduction of an (intrinsic) electromagnetic
field on $T'M$, with the $2n$-component real quantity
\begin{equation}
\label{electromagnetic covector potential}
i\,(\,[v^{(j)}]^*v^{(j)}_{x^\alpha}\,,\,[v^{(j)}]^*v^{(j)}_{\xi_\gamma}\,)
\end{equation}
playing the role of the electromagnetic covector potential. Our
quantity (\ref{electromagnetic covector potential}) is a 1-form on
$T'M$, rather than on $M$ itself as is the case in ``traditional''
electromagnetism. The above $\mathrm{U}(1)$ connection generates
curvature which is a 2-form on $T'M$, an analogue of the
electromagnetic tensor. Out of this curvature 2-form one can
construct, by contraction of indices, a real scalar. This scalar
curvature is the expression appearing in the RHS of
formula (\ref{subprincipal symbol of OI at time zero}).

Observe now that $\sum_jU^{(j)}(0)$ is the identity operator
on half-densities. The subprincipal symbol of the identity operator
is zero, so formula (\ref{subprincipal symbol of OI at time zero})
implies
\begin{equation}
\label{sum of curvatures is zero}
\sum_j\{[v^{(j)}]^*,v^{(j)}\}=0.
\end{equation}
One can check the identity (\ref{sum of curvatures is zero})
directly, without constructing the oscillatory integrals
$U^{(j)}(t)$: it follows from the fact that the $v^{(j)}(x,\xi)$ form an
orthonormal basis, see end of Section \ref{U(1) connection} for details.
We mentioned the identity
(\ref{sum of curvatures is zero})
in order to highlight, once again, the fact that the curvature effects
we have identified are specific to systems and do not have an
analogue in the scalar case.

\

\textbf{Results 2 and 3.}
We prove, under appropriate assumptions on Hamiltonian trajectories
(see Theorems~\ref{theorem spectral function unmollified two term}
and \ref{theorem counting function unmollified two term}),
asymptotic formulae
(\ref{two-term asymptotic formula for spectral function})
and
(\ref{two-term asymptotic formula for counting function})
with
\begin{equation}
\label{formula for a(x)}
a(x)=\sum_{j=1}^{m^+}
\ \int\limits_{h^{(j)}(x,\xi)<1}{d{\hskip-1pt\bar{}}\hskip1pt}\xi\,,
\end{equation}
\begin{multline}
\label{formula for b(x)}
b(x)=-n\sum_{j=1}^{m^+}
\ \int\limits_{h^{(j)}(x,\xi)<1}
\Bigl(
[v^{(j)}]^*A_\mathrm{sub}v^{(j)}
\\
-\frac i2
\{
[v^{(j)}]^*,A_1-h^{(j)},v^{(j)}
\}
+\frac i{n-1}h^{(j)}\{[v^{(j)}]^*,v^{(j)}\}
\Bigr)(x,\xi)\,
{d{\hskip-1pt\bar{}}\hskip1pt}\xi\,,
\end{multline}
and $a$ and $b$ expressed via the above densities
(\ref{formula for a(x)})
and
(\ref{formula for b(x)})
as
(\ref{a via a(x)})
and
(\ref{b via b(x)}).
In
(\ref{formula for a(x)})
and
(\ref{formula for b(x)})
\,${d{\hskip-1pt\bar{}}\hskip1pt}\xi$ is shorthand for
${d{\hskip-1pt\bar{}}\hskip1pt}\xi:=(2\pi)^{-n}\,d\xi
=(2\pi)^{-n}\,d\xi_1\ldots d\xi_n$,
and the Poisson bracket on matrix-functions
$\{\,\cdot\,,\,\cdot\,\}$
and its further generalisation
$\{\,\cdot\,,\,\cdot\,,\,\cdot\,\}$
are defined by formulae
(\ref{Poisson bracket on matrix-functions})
and
(\ref{generalised Poisson bracket on matrix-functions})
respectively.

To our knowledge, formula (\ref{formula for b(x)}) is a new result.
Note that in \cite{SafarovDSc} this formula
(more precisely, its integrated over $M$ version (\ref{b via b(x)}))
was written incorrectly,
without the curvature terms
$\,-\frac{ni}{n-1}\int h^{(j)}\{[v^{(j)}]^*,v^{(j)}\}$.
See also Section~\ref{Bibliographic review}
where we give a more detailed bibliographic review.

It is easy to see that the right-hand sides of
(\ref{formula for a(x)})
and
(\ref{formula for b(x)})
behave as densities under changes of local coordinates
on the manifold $M$ and that these expressions are invariant
under gauge transformations
(\ref{gauge transformation of the eigenvector}),
(\ref{phase appearing in gauge transformation})
of the eigenvectors of the principal symbol.
Moreover, the right-hand sides of
(\ref{formula for a(x)})
and
(\ref{formula for b(x)})
are unitarily invariant,
i.e.~invariant under the transformation of the operator
\begin{equation}
\label{unitary transformation of operator A}
A\mapsto RAR^*,
\end{equation}
where
\begin{equation}
\label{matrix appearing in unitary transformation of operator}
R:M\to\mathrm{U}(m)
\end{equation}
is an arbitrary smooth unitary matrix-function.
The fact that the RHS of
(\ref{formula for b(x)})
is unitarily invariant is non-trivial: the appropriate calculations
are presented in Section~\ref{U(m) invariance}.
The observation that without the curvature terms
$\,-\frac{ni}{n-1}\int h^{(j)}\{[v^{(j)}]^*,v^{(j)}\}$
(as in \cite{SafarovDSc}) the RHS of
(\ref{formula for b(x)}) is not unitarily invariant
was a major motivating factor in the writing of this paper.

\

Formula (\ref{formula for b(x)}) is the main result of our paper.
Note that even though the two-term asymptotic
expansion (\ref{two-term asymptotic formula for spectral function})
holds only under certain assumptions
on Hamiltonian trajectories (loops),
the second asymptotic coefficient (\ref{formula for b(x)})
is, in itself,
well-defined irrespective of how many loops we have. If one wishes
to reformulate the asymptotic expansion
(\ref{two-term asymptotic formula for spectral function})
in such a way that it remains valid without assumptions on the
number of loops, this can easily be achieved, say, by
taking a convolution with a function
from Schwartz space $\mathcal{S}(\mathbb{R})$:
see Theorem~\ref{theorem spectral function mollified}.

%

\section{Algorithm for the construction of the propagator}
\label{Algorithm for the construction of the wave group}

We construct the propagator as a sum of $m$ oscillatory integrals
(\ref{wave group as a sum of oscillatory integrals}) where each
integral is of the form
\begin{equation}
\label{algorithm equation 1}
U^{(j)}(t)
=
\int e^{i\varphi^{(j)}(t,x;y,\eta)}
\,u^{(j)}(t;y,\eta)
\,\varsigma^{(j)}(t,x;y,\eta)\,d_{\varphi^{(j)}}(t,x;y,\eta)\,
(\ \cdot\ )\,dy\,{d{\hskip-1pt\bar{}}\hskip1pt}\eta\,.
\end{equation}
Here we use notation from the book \cite{mybook}, only adapted to
systems. Namely, the expressions appearing in formula
(\ref{algorithm equation 1}) have the following meaning.

\begin{itemize}
\item
The function
$\varphi^{(j)}$ is a phase function, i.e.~a function
$\mathbb{R}\times M\times T'M\to\mathbb{C}$
positively homogeneous in $\eta$ of degree 1 and satisfying
the conditions
\begin{equation}
\label{algorithm equation 2}
\varphi^{(j)}(t,x;y,\eta)
=(x-x^{(j)}(t;y,\eta))^\alpha\,\xi^{(j)}_\alpha(t;y,\eta)
+O(|x-x^{(j)}(t;y,\eta)|^2),
\end{equation}
\begin{equation}
\label{algorithm equation 3}
\operatorname{Im}\varphi^{(j)}(t,x;y,\eta)\ge0,
\end{equation}
\begin{equation}
\label{algorithm equation 4}
\det\varphi^{(j)}_{x^\alpha\eta_\beta}(t,x^{(j)}(t;y,\eta);y,\eta)\ne0.
\end{equation}
Recall that according to Corollary 2.4.5 from
\cite{mybook} we are guaranteed to have
(\ref{algorithm equation 4}) if we choose a phase function
\begin{multline}
\label{algorithm equation 5}
\varphi^{(j)}(t,x;y,\eta)
=(x-x^{(j)}(t;y,\eta))^\alpha\,\xi^{(j)}_\alpha(t;y,\eta)
\\
+\frac12C^{(j)}_{\alpha\beta}(t;y,\eta)
\,(x-x^{(j)}(t;y,\eta))^\alpha\,(x-x^{(j)}(t;y,\eta))^\beta
\\
+O(|x-x^{(j)}(t;y,\eta)|^3)
\end{multline}
with complex-valued symmetric matrix-function $C^{(j)}_{\alpha\beta}$ satisfying the strict inequality
$\operatorname{Im}C^{(j)}>0$
(our original requirement (\ref{algorithm equation 3}) implies only the
non-strict inequality $\operatorname{Im}C^{(j)}\ge0$).
Note that even though the matrix-function $C^{(j)}_{\alpha\beta}$ is
not a tensor, the inequalities $\operatorname{Im}C^{(j)}\ge0$ and
$\operatorname{Im}C^{(j)}>0$ are invariant under transformations of
local coordinates $x$; see Remark 2.4.9 in \cite{mybook} for details.

\item
The quantity $u^{(j)}$ is the symbol of our oscillatory integral,
i.e.~a complex-valued $m\times m$ matrix-function
$\mathbb{R}\times T'M\to\mathbb{C}^{m^2}$
which admits the asymptotic expansion
(\ref{decomposition of symbol of OI into homogeneous components}).
The symbol is the unknown quantity in our construction.

\item
The quantity $d_{\varphi^{(j)}}$ is defined in accordance with
formula (2.2.4) from \cite{mybook} as
\begin{equation}
\label{algorithm equation 6}
d_{\varphi^{(j)}}(t,x;y,\eta)
:=({\det}^2\varphi^{(j)}_{x^\alpha\eta_\beta})^{1/4}
=|\det\varphi^{(j)}_{x^\alpha\eta_\beta}|^{1/2}
\,e^{\,i\arg({\det}^2\varphi^{(j)}_{x^\alpha\eta_\beta})/4}.
\end{equation}
Note that in view of (\ref{algorithm equation 4}) our $d_{\varphi^{(j)}}$
is well-defined and smooth for $x$ close to $x^{(j)}(t;y,\eta)$. It is known
\cite{mybook} that under coordinate transformations
$d_{\varphi^{(j)}}$ behaves as a half-density in $x$ and
as a half-density to the power $-1$ in $y$.

In formula (\ref{algorithm equation 6}) we wrote
$({\det}^2\varphi^{(j)}_{x^\alpha\eta_\beta})^{1/4}$
rather than
$(\det\varphi^{(j)}_{x^\alpha\eta_\beta})^{1/2}$
in order to make this expression truly invariant under coordinate transformations.
Recall that local coordinates $x$ and $y$ are chosen independently
and that $\eta$ is a covector based at the point $y$.
Consequently,
$\det\varphi^{(j)}_{x^\alpha\eta_\beta}$ changes sign under inversion
of one of the local coordinates $\,x^\alpha$, $\alpha=1,\ldots,n$,
or $\,y^\beta$, $\beta=1,\ldots,n$,
whereas ${\det}^2\varphi^{(j)}_{x^\alpha\eta_\beta}$ retains sign under
inversion.

The choice of (smooth) branch of $\arg({\det}^2\varphi^{(j)}_{x^\alpha\eta_\beta})$ is assumed
to be fixed. Thus, for a given phase function
$\varphi^{(j)}$
formula (\ref{algorithm equation 6}) defines the quantity
$d_{\varphi^{(j)}}$
uniquely up to a factor $e^{ik\pi/2}$, $k=0,1,2,3$.
Observe now that if we set $t=0$ and choose the same local coordinates
for $x$ and $y$, we get $\varphi^{(j)}_{x^\alpha\eta_\beta}(0,y;y,\eta)=I$.
This implies that we can fully specify the choice of branch of
$\arg({\det}^2\varphi^{(j)}_{x^\alpha\eta_\beta})$
by requiring that
$d_{\varphi^{(j)}}(0,y;y,\eta)=1$.

The purpose of the introduction of the factor $d_{\varphi^{(j)}}$
in (\ref{algorithm equation 1}) is twofold.
\begin{itemize}
\item[(a)]
It ensures that the symbol $u^{(j)}$ is a function on
$\mathbb{R}\times T'M$ in the full differential geometric sense of the word,
i.e.~that it is invariant under transformations of local coordinates $x$ and $y$.
\item[(b)]
It ensures that the principal symbol $u^{(j)}_0$ does not depend
on the choice of phase function $\varphi^{(j)}$.
See Remark 2.2.8 in \cite{mybook} for more details.
\end{itemize}

\item
The quantity $\varsigma^{(j)}$ is a smooth cut-off function
$\mathbb{R}\times M\times T'M\to\mathbb{R}$
satisfying the following conditions.
\begin{itemize}
\item[(a)]
$\varsigma^{(j)}(t,x;y,\eta)=0$ on the set
$\{(t,x;y,\eta):\ |h^{(j)}(y,\eta)|\le1/2\}$.
\item[(b)]
$\varsigma^{(j)}(t,x;y,\eta)=1$ on the intersection
of a small conic neighbourhood of the set
\begin{equation}
\label{algorithm equation 7.1}
\{(t,x;y,\eta):\ x=x^{(j)}(t;y,\eta)\}
\end{equation}
with the set $\{(t,x;y,\eta):\ |h^{(j)}(y,\eta)|\ge1\}$.
\item[(c)]
$\varsigma^{(j)}(t,x;y,\lambda\eta)=\varsigma^{(j)}(t,x;y,\eta)$
for $\,|h^{(j)}(y,\eta)|\ge1$, $\,\lambda\ge1$.
\end{itemize}

\item
It is known (see Section 2.3 in \cite{mybook} for details)
that Hamiltonian trajectories generated by a Hamiltonian
$h^{(j)}(x,\xi)$ positively homogeneous in $\xi$ of degree~1
satisfy the identity
\begin{equation}
\label{algorithm equation 7.1.5}
(x^{(j)}_\eta)^{\alpha\beta}\xi^{(j)}_\alpha=0,
\end{equation}
where $(x^{(j)}_\eta)^{\alpha\beta}:=\partial(x^{(j)})^\alpha/\partial\eta_\beta$.
Formulae (\ref{algorithm equation 2}) and (\ref{algorithm equation 7.1.5})
imply
\begin{equation}
\label{algorithm equation 7.2}
\varphi^{(j)}_\eta(t,x^{(j)}(t;y,\eta);y,\eta)=0.
\end{equation}
This allows us to apply the stationary phase method in the neighbourhood
of the set (\ref{algorithm equation 7.1}) and disregard what happens
away from it.
\end{itemize}

\

Our task now is to construct the symbols  $u^{(j)}_0(t;y,\eta)$, $j=1,\ldots,m$,
so that our oscillatory integrals $U^{(j)}(t)$, $j=1,\ldots,m$,
satisfy the dynamic equations
\begin{equation}
\label{algorithm equation 8}
(D_t+A(x,D_x))\,U^{(j)}(t)\overset{\operatorname{mod}C^\infty}=0
\end{equation}
and initial condition
\begin{equation}
\label{algorithm equation 9}
\sum_jU^{(j)}(0)\overset{\operatorname{mod}C^\infty}=I\,,
\end{equation}
where $I$ is the identity operator on half-densities;
compare with formulae
(\ref{dynamic equation most basic}),
(\ref{initial condition most basic})
and (\ref{wave group as a sum of oscillatory integrals}).
Note that the pseudodifferential operator $A$ in formula
(\ref{algorithm equation 8}) acts on the oscillatory integral
$U(t)$ in the variable $x$; say, if $A$ is a differential
operator this means that in order to evaluate $A\,U^{(j)}(t)$
one has to perform the appropriate
differentiations of the oscillatory integral
(\ref{algorithm equation 1})
in the variable $x$.
Following the conventions of Section 3.3 of \cite{mybook},
we emphasise the fact that the pseudodifferential operator $A$ in formula
(\ref{algorithm equation 8}) acts on the oscillatory integral
$U(t)$ in the variable $x$ by writing this pseudodifferential operator
as $A(x,D_x)$, where
$D_{x^\alpha}:=-i\partial/\partial x^\alpha$.

We examine first the dynamic equation (\ref{algorithm equation 8}).
We have
\[
(D_t+A(x,D_x))\,U^{(j)}(t)=F^{(j)}(t)\,,
\]
where $F^{(j)}(t)$ is the oscillatory integral
\[
F^{(j)}(t)
=
\int e^{i\varphi^{(j)}(t,x;y,\eta)}
\,f^{(j)}(t,x;y,\eta)
\,\varsigma^{(j)}(t,x;y,\eta)\,d_{\varphi^{(j)}}(t,x;y,\eta)\,
(\ \cdot\ )\,dy\,{d{\hskip-1pt\bar{}}\hskip1pt}\eta
\]
whose matrix-valued amplitude $f^{(j)}$ is given by the formula
\begin{equation}
\label{algorithm equation 12}
f^{(j)}=D_tu^{(j)}+
\bigl(
\varphi^{(j)}_t+(d_{\varphi^{(j)}})^{-1}(D_t d_{\varphi^{(j)}})+s^{(j)}
\bigr)
\,u^{(j)},
\end{equation}
where the matrix-function $s^{(j)}(t,x;y,\eta)$ is defined as
\begin{equation}
\label{algorithm equation 13}
s^{(j)}=e^{-i\varphi^{(j)}}(d_{\varphi^{(j)}})^{-1}\,A(x,D_x)\,(e^{i\varphi^{(j)}}d_{\varphi^{(j)}})\,.
\end{equation}

Theorem 18.1 from \cite{shubin} gives us the following explicit asymptotic
(in inverse powers of $\eta$) formula for the
matrix-function (\ref{algorithm equation 13}):
\begin{equation}
\label{algorithm equation 14}
s^{(j)}=(d_{\varphi^{(j)}})^{-1}\sum_{\bm\alpha}
\frac1{{\bm\alpha}!}
\,A^{({\bm\alpha})}(x,\varphi^{(j)}_x)\,(D_z^{\bm\alpha}\chi^{(j)})\bigr|_{z=x}\ ,
\end{equation}
where
\begin{equation}
\label{algorithm equation 15}
\chi^{(j)}(t,z,x;y,\eta)
=e^{i\psi^{(j)}(t,z,x;y,\eta)}d_{\varphi^{(j)}}(t,z;y,\eta),
\end{equation}
\begin{equation}
\label{algorithm equation 16}
\psi^{(j)}(t,z,x;y,\eta)
=\varphi^{(j)}(t,z;y,\eta)
-\varphi^{(j)}(t,x;y,\eta)
-\varphi^{(j)}_{x^\beta}(t,x;y,\eta)\,(z-x)^\beta.
\end{equation}
In formula (\ref{algorithm equation 14})
\begin{itemize}
\item
${\bm\alpha}:=(\alpha_1,\ldots,\alpha_n)$ is a multi-index
(note the bold font which we use to distinguish
multi-indices and individual indices),
${\bm\alpha}!:=\alpha_1!\cdots\alpha_n!\,$,
$D_z^{\bm\alpha}:=D_{z^1}^{\alpha_1}\cdots D_{z^n}^{\alpha_n}$,
$D_{z^\beta}:=-i\partial/\partial z^\beta$,
\item
$A(x,\xi)$ is the full symbol of the pseudodifferential operator $A$
written in local coordinates~$x$,
\item
$A^{({\bm\alpha})}(x,\xi):=\partial_\xi^{\bm\alpha}A(x,\xi)$,
$\partial_\xi^{\bm\alpha}:=\partial_{\xi_1}^{\alpha_1}\cdots\partial_{\xi_n}^{\alpha_n}$
and $\partial_{\xi_\beta}:=\partial/\partial\xi_\beta\,$.
\end{itemize}

When $|\eta|\to+\infty$
the matrix-valued amplitude $f^{(j)}(t,x;y,\eta)$ defined by formula
(\ref{algorithm equation 12}) admits an asymptotic expansion
\begin{equation}
\label{algorithm equation 17}
f^{(j)}(t,x;y,\eta)=f^{(j)}_1(t,x;y,\eta)+f^{(j)}_0(t,x;y,\eta)+f^{(j)}_{-1}(t,x;y,\eta)+\ldots
\end{equation}
into components positively homogeneous in $\eta$, with the subscript
indicating degree of homogeneity. Note the following differences between formulae
(\ref{decomposition of symbol of OI into homogeneous components})
and (\ref{algorithm equation 17}).
\begin{itemize}
\item
The leading term in
(\ref{algorithm equation 17})
has degree of homogeneity 1, rather than 0 as in
(\ref{decomposition of symbol of OI into homogeneous components}).
In fact, the leading term in
(\ref{algorithm equation 17})
can be easily written out explicitly
\begin{equation}
\label{algorithm equation 18}
f^{(j)}_1(t,x;y,\eta)=
(\varphi^{(j)}_t(t,x;y,\eta)+A_1(x,\varphi^{(j)}_x(t,x;y,\eta)))\,u^{(j)}_0(t;y,\eta)\,,
\end{equation}
where $A_1(x,\xi)$ is the (matrix-valued) principal symbol of the pseudodifferential
operator $A$.
\item
Unlike the symbol $u^{(j)}(t;y,\eta)$, the amplitude
$f^{(j)}(t,x;y,\eta)$ depends on $x$.
\end{itemize}

We now need to exclude the dependence on $x$ from the amplitude
$f^{(j)}(t,x;y,\eta)$. This can be done by means of the algorithm
described in subsection 2.7.3 of \cite{mybook}.
We outline this algorithm below.

Working in local coordinates, define the matrix-function
$\varphi^{(j)}_{x\eta}$ in accordance with
$(\varphi^{(j)}_{x\eta})_\alpha{}^\beta:=\varphi^{(j)}_{x^\alpha\eta_\beta}$
and then define its inverse $(\varphi^{(j)}_{x\eta})^{-1}$ from the identity
$(\varphi^{(j)})_\alpha{}^\beta[(\varphi^{(j)}_{x\eta})^{-1}]_\beta{}^\gamma:=\delta_\alpha{}^\gamma$.
Define the ``scalar'' first order linear differential operators
\begin{equation}
\label{algorithm equation 19}
L^{(j)}_\alpha:=[(\varphi^{(j)}_{x\eta})^{-1}]_\alpha{}^\beta\,(\partial/\partial x^\beta),
\qquad\alpha=1,\ldots,n.
\end{equation}
Note that the coefficients of these differential operators are functions of the position
variable $x$ and the dual variable $\xi$. It is known, see part 2 of Appendix E in \cite{mybook},
that the operators (\ref{algorithm equation 19}) commute:
$\ L^{(j)}_\alpha L^{(j)}_\beta=L^{(j)}_\beta L^{(j)}_\alpha$,
$\ \alpha,\beta=1,\ldots,n$.

Denote
$\ L^{(j)}_{\bm\alpha}:=(L^{(j)}_1)^{\alpha_1}\cdots(L^{(j)}_n)^{\alpha_n}$,
$\ (-\varphi^{(j)}_\eta)^{\bm\alpha}:=(-\varphi^{(j)}_{\eta_1})^{\alpha_1}\cdots(-\varphi^{(j)}_{\eta_n})^{\alpha_n}$,
and, given an $r\in\mathbb{N}$, define the ``scalar'' linear differential operator
\begin{equation}
\label{algorithm equation 21}
\mathfrak{P}^{(j)}_{-1,r}:=
i(d_{\varphi^{(j)}})^{-1}
\,
\frac\partial{\partial\eta_\beta}
\,d_{\varphi^{(j)}}
\left(1+
\sum_{1\le|{\bm\alpha}|\le2r-1}
\frac{(-\varphi^{(j)}_\eta)^{\bm\alpha}}{{\bm\alpha}!\,(|{\bm\alpha}|+1)}
\,L^{(j)}_{\bm\alpha}
\right)
L^{(j)}_\beta\,,
\end{equation}
where $|{\bm\alpha}|:=\alpha_1+\ldots+\alpha_n$ and the repeated index $\beta$ indicates
summation over $\beta=1,\ldots,n$.

Recall Definition 2.7.8 from \cite{mybook}:
the linear operator $L$ is said to be
positively homogeneous in $\eta$ of degree $p\in\mathbb{R}$
if for any $q\in\mathbb{R}$ and any function $f$
positively homogeneous in $\eta$ of degree $q$
the function $Lf$ is
positively homogeneous in $\eta$ of degree $p+q$.
It is easy to see that the operator (\ref{algorithm equation 21}) is
positively homogeneous in $\eta$ of degree $-1$
and the first subscript in $\mathfrak{P}^{(j)}_{-1,r}$ emphasises this fact.

Let $\mathfrak{S}^{(j)}_0$ be the (linear) operator of restriction to $x=x^{(j)}(t;y,\eta)$,
\begin{equation}
\label{algorithm equation 22}
\mathfrak{S}^{(j)}_0:=\left.(\,\cdot\,)\right|_{x=x^{(j)}(t;y,\eta)}\,,
\end{equation}
and let
\begin{equation}
\label{algorithm equation 23}
\mathfrak{S}^{(j)}_{-r}:=\mathfrak{S}^{(j)}_0(\mathfrak{P}^{(j)}_{-1,r})^r
\end{equation}
for $r=1,2,\ldots$. Observe that our linear operators
$\mathfrak{S}^{(j)}_{-r}$, $r=0,1,2,\ldots$, are
positively homogeneous in $\eta$ of degree $-r$.
This observation allows us to define the linear operator
\begin{equation}
\label{algorithm equation 24}
\mathfrak{S}^{(j)}:=\sum_{r=0}^{+\infty}\mathfrak{S}^{(j)}_{-r}\ ,
\end{equation}
where the series is understood as an asymptotic series in inverse powers of $\eta$.

According to subsection 2.7.3 of \cite{mybook},
the dynamic equation (\ref{algorithm equation 8}) can now be rewritten in the equivalent form
\begin{equation}
\label{algorithm equation 25}
\mathfrak{S}^{(j)}f^{(j)}=0\,,
\end{equation}
where the equality is understood in the asymptotic sense, as
an asymptotic expansion in inverse powers of $\eta$.
Recall that the matrix-valued amplitude $f^{(j)}(t,x;y,\eta)$
appearing in (\ref{algorithm equation 25}) is defined
by formulae (\ref{algorithm equation 12})--(\ref{algorithm equation 16}).

Substituting (\ref{algorithm equation 24}) and (\ref{algorithm equation 17})
into (\ref{algorithm equation 25}) we obtain a hierarchy of equations
\begin{equation}
\label{algorithm equation 26}
\mathfrak{S}^{(j)}_0f^{(j)}_1=0,
\end{equation}
\begin{equation}
\label{algorithm equation 27}
\mathfrak{S}^{(j)}_{-1}f^{(j)}_1+\mathfrak{S}^{(j)}_0f^{(j)}_0=0,
\end{equation}
\[
\mathfrak{S}^{(j)}_{-2}f^{(j)}_1+\mathfrak{S}^{(j)}_{-1}f^{(j)}_0+\mathfrak{S}^{(j)}_0f^{(j)}_{-1}=0,
\]
\[
\ldots
\]
positively homogeneous in $\eta$ of degree 1, 0, $-1$, $\ldots$.
These are the \emph{transport} equations for the determination of the unknown
homogeneous components $u^{(j)}_0(t;y,\eta)$, $u^{(j)}_{-1}(t;y,\eta)$, $u^{(j)}_{-2}(t;y,\eta)$, $\ldots$,
of the symbol of the oscillatory integral (\ref{algorithm equation 1}).

Let us now examine the initial condition (\ref{algorithm equation 9}).
Each operator $U^{(j)}(0)$ is a pseudodifferential operator, only
written in a slightly nonstandard form. The issues here are as follows.
\begin{itemize}
\item
We use the invariantly defined phase function
$
\varphi^{(j)}(0,x;y,\eta)
=(x-y)^\alpha\,\eta_\alpha
+O(|x-y|^2)
$
rather than the linear phase function $(x-y)^\alpha\,\eta_\alpha$
written in local coordinates.
\item
When defining the (full) symbol of the operator $U^{(j)}(t)$ we excluded the variable
$x$ from the amplitude rather than the variable $y$. Note that when dealing
with pseudodifferential operators it is customary to exclude the variable $y$
from the amplitude; exclusion of the variable $x$ gives the dual symbol of
a pseudodifferential operator, see subsection 2.1.3 in \cite{mybook}.
Thus, at $t=0$, our symbol $u^{(j)}(0;y,\eta)$ resembles
the dual symbol of a pseudodifferential operator rather
than the ``normal'' symbol.
\item
We have the extra factor $d_{\varphi^{(j)}}(0,x;y,\eta)$ in our representation
of the operator $U^{(j)}(0)$ as an oscillatory integral.
\end{itemize}

The (full) dual symbol
of the pseudodifferential operator $U^{(j)}(0)$
can be calculated in local coordinates in accordance with the following
formula which addresses the issues highlighted above:
\begin{equation}
\label{algorithm equation 30}
\sum_{\bm\alpha}
\frac{(-1)^{|{\bm\alpha}|}}{{\bm\alpha}!}\,
\bigl(
D_x^{\bm\alpha}\,\partial_\eta^{\bm\alpha}\,
u^{(j)}(0;y,\eta)\,
e^{i\omega^{(j)}(x;y,\eta)}\,d_{\varphi^{(j)}}(0,x;y,\eta)
\bigr)
\bigr|_{x=y}\ ,
\end{equation}
where
$\omega^{(j)}(x;y,\eta)=\varphi^{(j)}(0,x;y,\eta)-(x-y)^\beta\,\eta_\beta\,$.
Formula (\ref{algorithm equation 30})
is a version of the formula from subsection 2.1.3 of \cite{mybook}, only with
the extra factor $(-1)^{|{\bm\alpha}|}$. The latter is needed because we are writing
down the dual symbol of the pseudodifferential operator $U^{(j)}(0)$ (no dependence on $x$)
rather than its ``normal'' symbol (no dependence on $y$).

The initial condition (\ref{algorithm equation 9}) can now be rewritten in explicit form as
\begin{equation}
\label{algorithm equation 32}
\sum_j
\sum_{\bm\alpha}
\frac{(-1)^{|{\bm\alpha}|}}{{\bm\alpha}!}\,
\bigl(
D_x^{\bm\alpha}\,\partial_\eta^{\bm\alpha}\,
u^{(j)}(0;y,\eta)\,
e^{i\omega^{(j)}(x;y,\eta)}\,d_{\varphi^{(j)}}(0,x;y,\eta)
\bigr)
\bigr|_{x=y}=I\,,
\end{equation}
where $I$ is the $m\times m$ identity matrix.
Condition (\ref{algorithm equation 32})
can be decomposed into components positively homogeneous in $\eta$
of degree $0,-1,-2,\ldots$, giving us a hierarchy of initial conditions.
The leading (of degree of homogeneity 0) initial condition reads
\begin{equation}
\label{algorithm equation 33}
\sum_j
u^{(j)}_0(0;y,\eta)=I\,,
\end{equation}
whereas lower order initial conditions are more complicated
and depend on the choice of our phase functions $\varphi^{(j)}$.

\section{Leading transport equations}
\label{Leading transport equations}

Formulae
(\ref{algorithm equation 22}),
(\ref{algorithm equation 18}),
(\ref{algorithm equation 2}),
(\ref{Hamiltonian system of equations})
and the identity $\xi_\alpha h^{(j)}_{\xi_\alpha}(x,\xi)=h^{(j)}(x,\xi)$
(consequence of the fact that $h^{(j)}(x,\xi)$ is positively homogeneous in $\xi$ of degree~1)
give us the following explicit representation
for the leading transport equation
(\ref{algorithm equation 26}):
\begin{equation}
\label{Leading transport equations equation 1}
\!\!
\bigl[
A_1\bigl(x^{(j)}(t;y,\eta),\xi^{(j)}(t;y,\eta)\bigr)
-
h^{(j)}\bigl(x^{(j)}(t;y,\eta),\xi^{(j)}(t;y,\eta)\bigr)
\bigr]
\,u^{(j)}_0(t;y,\eta)=0.
\end{equation}
Here, of course,
$h^{(j)}\bigl(x^{(j)}(t;y,\eta),\xi^{(j)}(t;y,\eta)\bigr)=h^{(j)}(y,\eta)$.

Equation (\ref{Leading transport equations equation 1}) implies that
\begin{equation}
\label{Leading transport equations equation 2}
u^{(j)}_0(t;y,\eta)=v^{(j)}(x^{(j)}(t;y,\eta),\xi^{(j)}(t;y,\eta))
\,[w^{(j)}(t;y,\eta)]^T,
\end{equation}
where $v^{(j)}(z,\zeta)$ is the normalised eigenvector of the principal
symbol $A_1(z,\zeta)$ corresponding to the eigenvalue $h^{(j)}(z,\zeta)$
and $w^{(j)}:\mathbb{R}\times T'M\to\mathbb{C}^m$ is a column-function,
positively homogeneous in $\eta$ of degree 0, that remains to be found.
Formulae
(\ref{algorithm equation 33})
and
(\ref{Leading transport equations equation 2})
imply the following initial condition for the unknown column-function $w^{(j)}$:
\begin{equation}
\label{Leading transport equations equation 3}
w^{(j)}(0;y,\eta)=\overline{v^{(j)}(y,\eta)}.
\end{equation}

We now consider the
next transport equation in our hierarchy,
equation (\ref{algorithm equation 27}).
We will write down the two terms appearing in
(\ref{algorithm equation 27}) separately.

In view of formulae
(\ref{algorithm equation 18})
and
(\ref{algorithm equation 21})--(\ref{algorithm equation 23}),
the first term in (\ref{algorithm equation 27}) reads
\begin{multline}
\label{Leading transport equations equation 4}
\mathfrak{S}^{(j)}_{-1}f^{(j)}_1=
\\
i
\left.
\left[
(d_{\varphi^{(j)}})^{-1}
\frac\partial{\partial\eta_\beta}
d_{\varphi^{(j)}}
\left(1-
\frac12
\varphi^{(j)}_{\eta_\alpha}
L^{(j)}_\alpha
\right)
\left(
L^{(j)}_\beta
\bigl(\varphi^{(j)}_t+A_1(x,\varphi^{(j)}_x)\bigr)
\right)
u^{(j)}_0
\right]
\right|_{x=x^{(j)}}\,,
\end{multline}
where we dropped, for the sake of brevity,
the arguments $(t;y,\eta)$ in $u^{(j)}_0$ and $x^{(j)}$,
and the arguments $(t,x;y,\eta)$
in $\varphi^{(j)}_t$, $\varphi^{(j)}_x$, $\varphi^{(j)}_\eta$ and $d_{\varphi^{(j)}}\,$.
Recall that the differential operators $L^{(j)}_\alpha$ are defined in accordance with
formula (\ref{algorithm equation 19})
and the coefficients of these operators depend on $(t,x;y,\eta)$.

In view of formulae
(\ref{algorithm equation 12})--(\ref{algorithm equation 17})
and
(\ref{algorithm equation 22}),
the second term in (\ref{algorithm equation 27}) reads
\begin{multline}
\label{Leading transport equations equation 5}
\mathfrak{S}^{(j)}_0f^{(j)}_0=
D_tu^{(j)}_0
\\
+\left.\left[
(d_{\varphi^{(j)}})^{-1}
\left(D_t+(A_1)_{\xi_\alpha}D_{x^\alpha}\right)
d_{\varphi^{(j)}}
+A_0
-\frac i2(A_1)_{\xi_\alpha\xi_\beta}C^{(j)}_{\alpha\beta}
\right]
\right|_{x=x^{(j)}}u^{(j)}_0
\\
+\bigl[A_1-h^{(j)}\bigr]u^{(j)}_{-1}\,,
\end{multline}
where
\begin{equation}
\label{Leading transport equations equation 6}
C^{(j)}_{\alpha\beta}:=\left.\varphi^{(j)}_{x^\alpha x^\beta}\right|_{x=x^{(j)}}
\end{equation}
is the matrix-function from
(\ref{algorithm equation 5}).
In formulae
(\ref{Leading transport equations equation 5})
and
(\ref{Leading transport equations equation 6})
we dropped, for the sake of brevity,
the arguments $(t;y,\eta)$ in $u^{(j)}_0$, $u^{(j)}_{-1}$, $C^{(j)}_{\alpha\beta}$ and $x^{(j)}$,
the arguments
$(x^{(j)}(t;y,\eta),\xi^{(j)}(t;y,\eta))$
in $A_0$, $A_1$, $(A_1)_{\xi_\alpha}$, $(A_1)_{\xi_\alpha\xi_\beta}$ and $h^{(j)}$,
and the arguments $(t,x;y,\eta)$
in $d_{\varphi^{(j)}}$ and $\varphi^{(j)}_{x^\alpha x^\beta}\,$.

Looking at
(\ref{Leading transport equations equation 4})
and
(\ref{Leading transport equations equation 5})
we see that the transport equation (\ref{algorithm equation 27}) has a complicated
structure.
Hence, in this section we choose not to perform the analysis
of the full equation (\ref{algorithm equation 27})
and analyse only one particular subequation of this equation.
Namely, observe that equation (\ref{algorithm equation 27})
is equivalent to $m$ subequations
\begin{equation}
\label{Leading transport equations equation 7}
\bigl[v^{(j)}\bigr]^*
\,
\bigl[\mathfrak{S}^{(j)}_{-1}f^{(j)}_1+\mathfrak{S}^{(j)}_0f^{(j)}_0\bigr]
=0,
\end{equation}
\begin{equation}
\label{Leading transport equations equation 8}
\bigl[v^{(l)}\bigr]^*
\,
\bigl[\mathfrak{S}^{(j)}_{-1}f^{(j)}_1+\mathfrak{S}^{(j)}_0f^{(j)}_0\bigr]
=0,
\qquad l\ne j,
\end{equation}
where we dropped, for the sake of brevity, the arguments
$(x^{(j)}(t;y,\eta),\xi^{(j)}(t;y,\eta))$
in $\bigl[v^{(j)}\bigr]^*$ and $\bigl[v^{(l)}\bigr]^*$.
In the remainder of this section we analyse (sub)equation
(\ref{Leading transport equations equation 7}) only.

Equation (\ref{Leading transport equations equation 7})
is simpler than each of the $m-1$ equations
(\ref{Leading transport equations equation 8})
for the following two reasons.

\begin{itemize}

\item
Firstly, the term
$\bigl[A_1-h^{(j)}\bigr]u^{(j)}_{-1}$
from
(\ref{Leading transport equations equation 5})
vanishes after multiplication by
$\bigl[v^{(j)}\bigr]^*$
from the left.
Hence, equation
(\ref{Leading transport equations equation 7})
does not contain $u^{(j)}_{-1}$.

\item
Secondly, if we substitute
(\ref{Leading transport equations equation 2})
into
(\ref{Leading transport equations equation 7}),
then the term with
\[
\partial[d_{\varphi^{(j)}}w^{(j)}(t;y,\eta)]^T/\partial\eta_\beta
\]
vanishes.
This follows from the fact that the scalar function
\[
\bigl[v^{(j)}\bigr]^*
\bigl(\varphi^{(j)}_t+A_1(x,\varphi^{(j)}_x)\bigr)
v^{(j)}
\]
has a second order zero, in the variable $x$, at $x=x^{(j)}(t;y,\eta)$.
Indeed, we have
\begin{multline*}
\left.
\left[
\frac\partial{\partial x^\alpha}
\bigl[v^{(j)}\bigr]^*
\bigl(\varphi^{(j)}_t+A_1(x,\varphi^{(j)}_x)\bigr)
v^{(j)}
\right]
\right|_{x=x^{(j)}}
\\
=
\bigl[v^{(j)}\bigr]^*
\left.
\left[
\bigl(\varphi^{(j)}_t+A_1(x,\varphi^{(j)}_x)\bigr)_{x^\alpha}
\right]
\right|_{x=x^{(j)}}
v^{(j)}
\\
=
\bigl[v^{(j)}\bigr]^*
\bigl(
-h^{(j)}_{x^\alpha}-C^{(j)}_{\alpha\beta}h^{(j)}_{\xi_\beta}
+(A_1)_{x^\alpha}+C^{(j)}_{\alpha\beta}(A_1)_{\xi_\beta}
\bigr)
v^{(j)}
\\
=
\bigl[v^{(j)}\bigr]^*(A_1)_{x^\alpha}v^{(j)}-h^{(j)}_{x^\alpha}
+
C^{(j)}_{\alpha\beta}
\bigl(
\bigl[v^{(j)}\bigr]^*(A_1)_{\xi_\beta}v^{(j)}-h^{(j)}_{\xi_\beta}
\bigr)
=0\,,
\end{multline*}
where in the last two lines we dropped,
for the sake of brevity,
the arguments
$(x^{(j)}(t;y,\eta),\xi^{(j)}(t;y,\eta))$
in $(A_1)_{x^\alpha}$, $(A_1)_{\xi_\beta}$,
$h^{(j)}_{x^\alpha}$, $h^{(j)}_{\xi_\beta}$,
and the argument $(t;y,\eta)$ in
$C^{(j)}_{\alpha\beta}$
(the latter is the
matrix-function
from
formulae (\ref{algorithm equation 5}) and (\ref{Leading transport equations equation 6})).
Throughout the above argument we used the fact that our
$\bigl[v^{(j)}\bigr]^*$ and $v^{(j)}$ do not depend on $x$:
their argument is
$(x^{(j)}(t;y,\eta),\xi^{(j)}(t;y,\eta))$.

\end{itemize}

Substituting
(\ref{Leading transport equations equation 4}),
(\ref{Leading transport equations equation 5})
and
(\ref{Leading transport equations equation 2})
into
(\ref{Leading transport equations equation 7})
we get
\begin{equation}
\label{Leading transport equations equation 9}
(D_t+p^{(j)}(t;y,\eta))\,[w^{(j)}(t;y,\eta)]^T=0\,,
\end{equation}
where
\begin{multline}
\label{Leading transport equations equation 10}
p^{(j)}=
i
\left.
[v^{(j)}]^*
\left[
\frac\partial{\partial\eta_\beta}
\left(1-
\frac12
\varphi^{(j)}_{\eta_\alpha}
L^{(j)}_\alpha
\right)
\left(
L^{(j)}_\beta
\bigl(\varphi^{(j)}_t+A_1(x,\varphi^{(j)}_x)\bigr)
\right)
v^{(j)}
\right]
\right|_{x=x^{(j)}}
\\
-i[v^{(j)}]^*\{v^{(j)},h^{(j)}\}
+\left.\left[
(d_{\varphi^{(j)}})^{-1}
\left(D_t+h^{(j)}_{\xi_\alpha}D_{x^\alpha}\right)
d_{\varphi^{(j)}}
\right]
\right|_{x=x^{(j)}}
\\
+[v^{(j)}]^*
\left(
A_0
-\frac i2(A_1)_{\xi_\alpha\xi_\beta}C^{(j)}_{\alpha\beta}
\right)
v^{(j)}.
\end{multline}

Note that the ordinary differential operator in the LHS of formula
(\ref{Leading transport equations equation 9}) is a scalar one, i.e.
it does not mix up the different components of the column-function
$w^{(j)}(t;y,\eta)$. The solution of the ordinary
differential equation
(\ref{Leading transport equations equation 9})
subject to the initial condition
(\ref{Leading transport equations equation 3})
is
\begin{equation}
\label{Leading transport equations equation 11}
w^{(j)}(t;y,\eta)=\overline{v^{(j)}(y,\eta)}
\exp\left(-i\int_0^tp^{(j)}(\tau;y,\eta)\,d\tau\right).
\end{equation}
Comparing formulae
(\ref{Leading transport equations equation 2}),
(\ref{Leading transport equations equation 11})
with formula
(\ref{formula for principal symbol of oscillatory integral})
we see that in order to prove the latter we need only to establish the
scalar identity
\begin{equation}
\label{Leading transport equations equation 12}
p^{(j)}(t;y,\eta)=q^{(j)}(x^{(j)}(t;y,\eta),\xi^{(j)}(t;y,\eta))\,,
\end{equation}
where $q^{(j)}$ is the function
(\ref{phase appearing in principal symbol}).
In view of the definitions of the quantities
$p^{(j)}$ and $q^{(j)}$,
see formulae
(\ref{Leading transport equations equation 10})
and
(\ref{phase appearing in principal symbol}),
and the definition of the subprincipal symbol
(\ref{definition of subprincipal symbol}),
proving the identity
(\ref{Leading transport equations equation 12})
reduces to proving the identity
\begin{multline}
\label{Leading transport equations equation 13}
\{
[v^{(j)}]^*,A_1-h^{(j)},v^{(j)}
\}
(x^{(j)},\xi^{(j)})
=
\\
-2
\left.
[v^{(j)}(x^{(j)},\xi^{(j)})]^*
\left[
\frac\partial{\partial\eta_\beta}
\left(1-
\frac12
\varphi^{(j)}_{\eta_\alpha}
L^{(j)}_\alpha
\right)
\left(
L^{(j)}_\beta
\bigl(\varphi^{(j)}_t+A_1(x,\varphi^{(j)}_x)\bigr)
\right)
v^{(j)}(x^{(j)},\xi^{(j)})
\right]
\right|_{x=x^{(j)}}
\\
+2\left.\left[
(d_{\varphi^{(j)}})^{-1}
\left(\partial_t+h^{(j)}_{\xi_\alpha}\partial_{x^\alpha}\right)
d_{\varphi^{(j)}}
\right]
\right|_{x=x^{(j)}}
\\
+[v^{(j)}(x^{(j)},\xi^{(j)})]^*
\left(
(A_1)_{x^\alpha\xi_\alpha}+(A_1)_{\xi_\alpha\xi_\beta}C^{(j)}_{\alpha\beta}
\right)
v^{(j)}(x^{(j)},\xi^{(j)}).
\end{multline}
Note that the expressions in the LHS and RHS of
(\ref{Leading transport equations equation 13}) have different
structure. The LHS of
(\ref{Leading transport equations equation 13})
is the generalised Poisson bracket
$\{
[v^{(j)}]^*,A_1-h^{(j)},v^{(j)}
\}$,
see (\ref{generalised Poisson bracket on matrix-functions}),
evaluated at $z=x^{(j)}(t;y,\eta)$, $\zeta=\xi^{(j)}(t;y,\eta)$,
whereas the RHS of
(\ref{Leading transport equations equation 13})
involves partial derivatives (in $\eta$) of
$v^{(j)}(x^{(j)}(t;y,\eta),\xi^{(j)}(t;y,\eta))$
(Chain Rule).
In writing
(\ref{Leading transport equations equation 13}) we also
dropped, for the sake of brevity,
the arguments $(t,x;y,\eta)$ in
$\varphi^{(j)}_t$, $\varphi^{(j)}_x$,
$\varphi^{(j)}_\eta$, $d_{\varphi^{(j)}}\,$
and the coefficients of the differential operators
$L^{(j)}_\alpha$ and $L^{(j)}_\beta$,
the arguments $(x^{(j)},\xi^{(j)})$
in $h^{(j)}_{\xi_\alpha}$,
$(A_1)_{x^\alpha\xi_\alpha}$ and $(A_1)_{\xi_\alpha\xi_\beta}$,
and the arguments $(t;y,\eta)$
in $x^{(j)}$, $\xi^{(j)}$ and $C^{(j)}_{\alpha\beta}$.

Before performing the calculations that will establish the identity
(\ref{Leading transport equations equation 13}) we make several
observations that will allow us to simplify these calculations
considerably.

Firstly, our function $p^{(j)}(t;y,\eta)$ does not depend on the choice
of the phase function $\varphi^{(j)}(t,x;y,\eta)$. Indeed, if
$p^{(j)}(t;y,\eta)$ did depend on the choice of phase function, then,
in view of formulae
(\ref{Leading transport equations equation 2})
and
(\ref{Leading transport equations equation 11})
the principal symbol of our oscillatory integral $U^{(j)}(t)$ would depend
on the choice of phase function, which would contradict Theorem
2.7.11 from \cite{mybook}. Here we use the fact that operators
$U^{(j)}(t)$ with different $j$ cannot compensate each other to give
an integral operator whose integral kernel is infinitely
smooth in $t$, $x$ and $y$ because all our $U^{(j)}(t)$ oscillate in
$t$ in a different way:
$\varphi^{(j)}_t(t,x^{(j)}(t;y,\eta);y,\eta)=-h^{(j)}(y,\eta)$
and we assumed the eigenvalues $h^{(j)}(y,\eta)$ of our principal
symbol $A_1(y,\eta)$ to be simple.

Secondly, the arguments (free variables) in
(\ref{Leading transport equations equation 13}) are
$(t;y,\eta)$. We fix an arbitrary point
$(\tilde t;\tilde y,\tilde\eta)\in\mathbb{R}\times T'M$
and prove formula
(\ref{Leading transport equations equation 13})
at this point.
Put
$(\xi^{(j)}_\eta)_\alpha{}^\beta:=\partial(\xi^{(j)})_\alpha/\partial\eta_\beta$.
According to Lemma 2.3.2 from \cite{mybook}
there exists a local coordinate system $x$ such that
$\det(\xi^{(j)}_\eta)_\alpha{}^\beta\ne0$.
This opens the way to the use of the linear phase function
\begin{equation}
\label{Leading transport equations equation 14}
\varphi^{(j)}(t,x;y,\eta)
=(x-x^{(j)}(t;y,\eta))^\alpha\,\xi^{(j)}_\alpha(t;y,\eta)
\end{equation}
which will simplify calculations to a great extent.
Moreover, we can choose a local coordinate system $y$ such that
\begin{equation}
\label{Leading transport equations equation 15}
(\xi^{(j)}_\eta)_\alpha{}^\beta(\tilde t;\tilde y,\tilde\eta)=\delta_\alpha{}^\beta
\end{equation}
which will simplify calculations even further.

The calculations we are about to perform will make use of the symmetry
\begin{equation}
\label{Leading transport equations equation 16}
(x^{(j)}_\eta)^{\gamma\alpha}(\xi^{(j)}_\eta)_\gamma{}^\beta
=
(x^{(j)}_\eta)^{\gamma\beta}(\xi^{(j)}_\eta)_\gamma{}^\alpha
\end{equation}
which is an immediate
consequence of formula (\ref{algorithm equation 7.1.5}).
Formula (\ref{Leading transport equations equation 16}) appears
as formula (2.3.3) in \cite{mybook} and the accompanying text
explains its geometric meaning. Note that at the point
$(\tilde t;\tilde y,\tilde\eta)$ formula
(\ref{Leading transport equations equation 16}) takes the
especially simple form
\begin{equation}
\label{Leading transport equations equation 17}
(x^{(j)}_\eta)^{\alpha\beta}(\tilde t;\tilde y,\tilde\eta)
=
(x^{(j)}_\eta)^{\beta\alpha}(\tilde t;\tilde y,\tilde\eta).
\end{equation}

Our calculations will also involve the quantity
$\varphi^{(j)}_{\eta_\alpha\eta_\beta}(\tilde t,\tilde x;\tilde y,\tilde\eta)$
where $\tilde x:=x^{(j)}(\tilde t;\tilde y,\tilde\eta)$.
Formulae
(\ref{Leading transport equations equation 14}),
(\ref{algorithm equation 7.1.5}),
(\ref{Leading transport equations equation 15})
and
(\ref{Leading transport equations equation 17})
imply
\begin{equation}
\label{Leading transport equations equation 18}
\varphi^{(j)}_{\eta_\alpha\eta_\beta}(\tilde t,\tilde x;\tilde y,\tilde\eta)
=
-(x^{(j)}_\eta)^{\alpha\beta}(\tilde t;\tilde y,\tilde\eta).
\end{equation}

Further on we denote $\tilde\xi:=\xi^{(j)}(\tilde t;\tilde y,\tilde\eta)$.

With account of all the simplifications listed above, we can rewrite
formula
(\ref{Leading transport equations equation 13}),
which is the identity that we are proving, as
\begin{multline}
\label{Leading transport equations equation 19}
\{
[v^{(j)}]^*,A_1-h^{(j)},v^{(j)}
\}
(\tilde x,\tilde\xi)
=
\\
\!\!\!\!\!\!\!\!\!\!\!\!\!\!\!\!\!\!\!\!\!\!\!\!\!\!\!\!
-2
[\tilde v^{(j)}]^*
\Bigl[
\frac{\partial^2}{\partial x^\alpha\partial\eta_\alpha}
\bigl(A_1(x,\xi^{(j)})-h^{(j)}(\tilde y,\eta)
\\
\qquad\qquad\qquad\qquad\qquad\qquad
-(x-x^{(j)})^\gamma h^{(j)}_{x^\gamma}(x^{(j)},\xi^{(j)})\bigr)
\,v^{(j)}(x^{(j)},\xi^{(j)})
\Bigr]
\Bigr|_{(x,\eta)=(\tilde x,\tilde\eta)}
\\
\!\!\!\!\!\!\!\!\!\!\!\!\!\!\!\!\!\!\!\!\!\!\!\!\!\!\!\!\!\!\!\!\!\!\!\!\!\!\!\!\!\!\!\!\!\!\!
-(\tilde x^{(j)}_\eta)^{\alpha\beta}\,
[\tilde v^{(j)}]^*
\Bigl[
\frac{\partial^2}{\partial x^\alpha\partial x^\beta}
\bigl(A_1(x,\xi^{(j)})-h^{(j)}(\tilde y,\eta)
\\
\qquad\qquad\qquad\qquad\qquad\qquad
-(x-x^{(j)})^\gamma h^{(j)}_{x^\gamma}(x^{(j)},\xi^{(j)})\bigr)
\,v^{(j)}(x^{(j)},\xi^{(j)})
\Bigr]
\Bigr|_{(x,\eta)=(\tilde x,\tilde\eta)}
\\
+[\tilde v^{(j)}]^*
(\tilde A_1)_{x^\alpha\xi_\alpha}
\tilde v^{(j)}
-\tilde h^{(j)}_{x^\alpha\xi_\alpha}
-\tilde h^{(j)}_{x^\alpha x^\beta}(\tilde
x^{(j)}_\eta)^{\alpha\beta}\,,
\qquad
\end{multline}

\noindent
where
$\tilde v^{(j)}=v^{(j)}(\tilde x,\tilde\xi)$,
$\tilde x^{(j)}_\eta=x^{(j)}_\eta(\tilde t;\tilde y,\tilde\eta)$,
$(\tilde A_1)_{x^\alpha\xi_\alpha}=(A_1)_{x^\alpha\xi_\alpha}(\tilde x,\tilde\xi)$,
$\tilde h^{(j)}_{x^\alpha\xi_\alpha}=h^{(j)}_{x^\alpha\xi_\alpha}(\tilde x,\tilde\xi)$,
$\tilde h^{(j)}_{x^\alpha x^\beta}=h^{(j)}_{x^\alpha x^\beta}(\tilde x,\tilde\xi)$,
$x^{(j)}=x^{(j)}(\tilde t;\tilde y,\eta)$
and
$\xi^{(j)}=\xi^{(j)}(\tilde t;\tilde y,\eta)$.

Note that the last two terms in the RHS of
(\ref{Leading transport equations equation 19})
originate from the term with
$d_{\varphi^{(j)}}$ in
(\ref{Leading transport equations equation 13}):
we used the fact that $d_{\varphi^{(j)}}$ does not depend on
$x$ and that
\begin{equation}
\label{Leading transport equations equation 20}
\left.
\left[
(d_{\varphi^{(j)}})^{-1}
\partial_t
d_{\varphi^{(j)}}
\right]
\right|_{(t,x;y,\eta)=(\tilde t,\tilde x;\tilde y,\tilde\eta)}
=-\frac12
\bigl(
\tilde h^{(j)}_{x^\alpha\xi_\alpha}
+\tilde h^{(j)}_{x^\alpha x^\beta}(\tilde x^{(j)}_\eta)^{\alpha\beta}
\bigr).
\end{equation}
Formula (\ref{Leading transport equations equation 20})
is a special case of formula (3.3.21) from \cite{mybook}.

Note also that the term $-h^{(j)}(\tilde y,\eta)$
appearing (twice) in the RHS of
(\ref{Leading transport equations equation 19})
will vanish after being acted upon with
the differential operators
$\frac{\partial^2}{\partial x^\alpha\partial\eta_\alpha}$
and
$\frac{\partial^2}{\partial x^\alpha\partial x^\beta}$
because it does not depend on $x$.

We have
\begin{multline}
\label{Leading transport equations equation 21}
[\tilde v^{(j)}]^*
\left.
\left[
\frac{\partial^2}{\partial x^\alpha\partial\eta_\alpha}
\bigl(A_1(x,\xi^{(j)})-(x-x^{(j)})^\gamma h^{(j)}_{x^\gamma}(x^{(j)},\xi^{(j)})\bigr)
\,v^{(j)}(x^{(j)},\xi^{(j)})
\right]
\right|_{(x,\eta)=(\tilde x,\tilde\eta)}
\\
=
[\tilde v^{(j)}]^*
(\tilde A_1)_{x^\alpha\xi_\alpha}
\tilde v^{(j)}
-\tilde h^{(j)}_{x^\alpha\xi_\alpha}
-\tilde h^{(j)}_{x^\alpha x^\beta}(\tilde x^{(j)}_\eta)^{\alpha\beta}
\\
+
[\tilde v^{(j)}]^*
\bigl(
(\tilde A_1)_{x^\alpha}
-\tilde h^{(j)}_{x^\alpha}
\bigr)
\bigl(
\tilde v^{(j)}_{\xi_\alpha}
+\tilde v^{(j)}_{x^\beta}(\tilde x^{(j)}_\eta)^{\alpha\beta}
\bigr),
\end{multline}
\begin{multline}
\label{Leading transport equations equation 22}
[\tilde v^{(j)}]^*
\left.
\left[
\frac{\partial^2}{\partial x^\alpha\partial x^\beta}
\bigl(A_1(x,\xi^{(j)})-(x-x^{(j)})^\gamma h^{(j)}_{x^\gamma}(x^{(j)},\xi^{(j)})\bigr)
\,v^{(j)}(x^{(j)},\xi^{(j)})
\right]
\right|_{(x,\eta)=(\tilde x,\tilde\eta)}
\\
=
[\tilde v^{(j)}]^*
(\tilde A_1)_{x^\alpha x^\beta}
\tilde v^{(j)}\,,
\end{multline}
where
$(\tilde A_1)_{x^\alpha}=(A_1)_{x^\alpha}(\tilde x,\tilde\xi)$,
$\tilde h^{(j)}_{x^\alpha}=h^{(j)}_{x^\alpha}(\tilde x,\tilde\xi)$,
$\tilde v^{(j)}_{\xi_\alpha}=v^{(j)}_{\xi_\alpha}(\tilde x,\tilde\xi)$
and
$\tilde v^{(j)}_{x^\beta}=v^{(j)}_{x^\beta}(\tilde x,\tilde\xi)$.
We also have
\begin{multline}
\label{Leading transport equations equation 23}
[\tilde v^{(j)}]^*
\bigl(
(\tilde A_1)_{x^\alpha}
-\tilde h^{(j)}_{x^\alpha}
\bigr)
\tilde v^{(j)}_{x^\beta}
+
[\tilde v^{(j)}]^*
\bigl(
(\tilde A_1)_{x^\beta}
-\tilde h^{(j)}_{x^\beta}
\bigr)
\tilde v^{(j)}_{x^\alpha}
\\
=
\tilde h^{(j)}_{x^\alpha x^\beta}
-
[\tilde v^{(j)}]^*
(\tilde A_1)_{x^\alpha x^\beta}
\tilde v^{(j)}.
\end{multline}
Using formulae
(\ref{Leading transport equations equation 23})
and
(\ref{Leading transport equations equation 17})
we can rewrite formula
(\ref{Leading transport equations equation 21})
as
\begin{multline}
\label{Leading transport equations equation 24}
[\tilde v^{(j)}]^*
\left.
\left[
\frac{\partial^2}{\partial x^\alpha\partial\eta_\alpha}
\bigl(A_1(x,\xi^{(j)})-(x-x^{(j)})^\gamma h^{(j)}_{x^\gamma}(x^{(j)},\xi^{(j)})\bigr)
\,v^{(j)}(x^{(j)},\xi^{(j)})
\right]
\right|_{(x,\eta)=(\tilde x,\tilde\eta)}
\\
=
[\tilde v^{(j)}]^*
(\tilde A_1)_{x^\alpha\xi_\alpha}
\tilde v^{(j)}
-\tilde h^{(j)}_{x^\alpha\xi_\alpha}
+
[\tilde v^{(j)}]^*
\bigl(
(\tilde A_1)_{x^\alpha}
-\tilde h^{(j)}_{x^\alpha}
\bigr)
\tilde v^{(j)}_{\xi_\alpha}
\\
-\frac12
\bigl(
[\tilde v^{(j)}]^*
(\tilde A_1)_{x^\alpha x^\beta}
\tilde v^{(j)}
+
\tilde h^{(j)}_{x^\alpha x^\beta}
\bigr)
(\tilde x^{(j)}_\eta)^{\alpha\beta}.
\end{multline}
Substituting
(\ref{Leading transport equations equation 24})
and
(\ref{Leading transport equations equation 22})
into
(\ref{Leading transport equations equation 19})
we see that all the terms with $(\tilde x^{(j)}_\eta)^{\alpha\beta}$ cancel out
and we get
\begin{multline}
\label{Leading transport equations equation 25}
\{
[v^{(j)}]^*,A_1-h^{(j)},v^{(j)}
\}
(\tilde x,\tilde\xi)
=
\\
-[\tilde v^{(j)}]^*
\bigl(
(\tilde A_1)_{x^\alpha\xi_\alpha}
-
\tilde h^{(j)}_{x^\alpha\xi_\alpha}
\bigr)
\tilde v^{(j)}
-2
[\tilde v^{(j)}]^*
\bigl(
(\tilde A_1)_{x^\alpha}
-\tilde h^{(j)}_{x^\alpha}
\bigr)
\tilde v^{(j)}_{\xi_\alpha}.
\end{multline}
Thus, the proof of the identity
(\ref{Leading transport equations equation 13})
has been reduced to the proof of the identity~(\ref{Leading transport equations equation 25}).

Observe now that formula
(\ref{Leading transport equations equation 25})
no longer has Hamiltonian trajectories present in it.
This means that we can drop all the tildes and
rewrite
(\ref{Leading transport equations equation 25})
as
\begin{multline}
\label{Leading transport equations equation 26}
\{
[v^{(j)}]^*,A_1-h^{(j)},v^{(j)}
\}
=
\\
-[v^{(j)}]^*
\bigl(
A_1
-h^{(j)}
\bigr)_{x^\alpha\xi_\alpha}
v^{(j)}
-2
[v^{(j)}]^*
\bigl(
A_1
-h^{(j)}
\bigr)_{x^\alpha}
v^{(j)}_{\xi_\alpha}\,,
\end{multline}
where the arguments are $(x,\xi)$.
We no longer need to restrict our consideration to the particular
point $(x,\xi)=(\tilde x,\tilde\xi)$:
if we prove
(\ref{Leading transport equations equation 26})
for an arbitrary $(x,\xi)\in T'M$
we will prove it for a particular
$(\tilde x,\tilde\xi)\in T'M$.

The proof of the identity
(\ref{Leading transport equations equation 26})
is straightforward. We note that
\begin{multline}
\label{Leading transport equations equation 27}
[v^{(j)}]^*
(A_1-h^{(j)})_{x^\alpha\xi_\alpha}
v^{(j)}=
\\
-
[v^{(j)}]^*
(A_1-h^{(j)})_{x^\alpha}
v^{(j)}_{\xi_\alpha}
-
[v^{(j)}]^*
(A_1-h^{(j)})_{\xi_\alpha}
v^{(j)}_{x^\alpha}
\end{multline}
and substituting
(\ref{Leading transport equations equation 27})
into
(\ref{Leading transport equations equation 26})
reduce the latter to the form
\begin{multline}
\label{Leading transport equations equation 28}
\{
[v^{(j)}]^*,A_1-h^{(j)},v^{(j)}
\}
=
\\
[v^{(j)}]^*
\bigl(
A_1
-h^{(j)}
\bigr)_{\xi_\alpha}
v^{(j)}_{x^\alpha}
-
[v^{(j)}]^*
\bigl(
A_1
-h^{(j)}
\bigr)_{x^\alpha}
v^{(j)}_{\xi_\alpha}.
\end{multline}
But
\begin{equation}
\label{Leading transport equations equation 29}
[v^{(j)}]^*
\bigl(
A_1
-h^{(j)}
\bigr)_{x^\alpha}
=
-
[v^{(j)}_{x^\alpha}]^*
\bigl(
A_1
-h^{(j)}
\bigr),
\end{equation}
\begin{equation}
\label{Leading transport equations equation 30}
[v^{(j)}]^*
\bigl(
A_1
-h^{(j)}
\bigr)_{\xi_\alpha}
=
-
[v^{(j)}_{\xi_\alpha}]^*
\bigl(
A_1
-h^{(j)}
\bigr).
\end{equation}
Substituting
(\ref{Leading transport equations equation 29})
and
(\ref{Leading transport equations equation 30})
into
(\ref{Leading transport equations equation 28})
we get
\[
\{
[v^{(j)}]^*,A_1-h^{(j)},v^{(j)}
\}
=
[v^{(j)}_{x^\alpha}]^*
\bigl(
A_1
-h^{(j)}
\bigr)
v^{(j)}_{\xi_\alpha}
-
[v^{(j)}_{\xi_\alpha}]^*
\bigl(
A_1
-h^{(j)}
\bigr)
v^{(j)}_{x^\alpha}
\]
which agrees with the definition of the generalised Poisson bracket
(\ref{generalised Poisson bracket on matrix-functions}).

\section{Proof of formula (\ref{subprincipal symbol of OI at time zero})}
\label{Proof of formula}

In this section we prove formula
(\ref{subprincipal symbol of OI at time zero}).
Our approach is as follows.

We write down explicitly the transport equations
(\ref{Leading transport equations equation 8}) at $t=0$,
i.e.
\begin{equation}
\label{Proof of formula equation 1}
\bigl[v^{(l)}\bigr]^*
\,
\left.
\bigl[\mathfrak{S}^{(j)}_{-1}f^{(j)}_1+\mathfrak{S}^{(j)}_0f^{(j)}_0\bigr]
\right|_{t=0}
=0,
\qquad l\ne j.
\end{equation}
We use the same local coordinates for $x$ and $y$ and we
assume all our phase functions to be linear, i.e.~we assume
that for each $j$ we have
(\ref{Leading transport equations equation 14}).
Using linear phase functions is justified for small $t$ because we
have
$(\xi^{(j)}_\eta)_\alpha{}^\beta(0;y,\eta)=\delta_\alpha{}^\beta$
and, hence, $\det\varphi^{(j)}_{x^\alpha\eta_\beta}(t,x;y,\eta)\ne0$
for small $t$. Writing down equations
(\ref{Proof of formula equation 1}) for linear phase functions
is much easier than for general phase functions
(\ref{algorithm equation 2}).

Using linear phase functions has the additional
advantage that the initial condition
(\ref{algorithm equation 32}) simplifies and reads now
$\sum_ju^{(j)}(0;y,\eta)=I$.
In view of
(\ref{decomposition of symbol of OI into homogeneous components}),
this implies, in particular, that
\begin{equation}
\label{Proof of formula equation 2}
\sum_j
u^{(j)}_{-1}(0)=0.
\end{equation}
Here and further on in this section we drop,
for the sake of brevity, the arguments
$(y,\eta)$ in $u^{(j)}_{-1}$.

Of course, the formula we are proving,
formula (\ref{subprincipal symbol of OI at time zero}),
does not depend
on our choice of phase functions. It is just easier to carry
out calculations for linear phase functions.

We will show that
(\ref{Proof of formula equation 1})
is a system of complex linear algebraic equations for the unknowns
$u^{(j)}_{-1}(0)$. The total number of equations
(\ref{Proof of formula equation 1}) is $m^2-m$. However, for each
$j$ and $l$
the LHS of (\ref{Proof of formula equation 1}) is a row of $m$
elements, so (\ref{Proof of formula equation 1}) is, effectively,
a system of $m(m^2-m)$ scalar equations.

Equation
(\ref{Proof of formula equation 2})
is a single matrix equation, so it is,
effectively,
a system of $m^2$ scalar equations.

Consequently, the system
(\ref{Proof of formula equation 1}),
(\ref{Proof of formula equation 2}) is, effectively,
a system of $m^3$ scalar equations.
This is exactly the number of unknown scalar elements
in the $m$ matrices $u^{(j)}_{-1}(0)$.

In the remainder of this section we write down explicitly
the LHS of (\ref{Proof of formula equation 1})
and solve the linear algebraic system
(\ref{Proof of formula equation 1}),
(\ref{Proof of formula equation 2})
for the unknowns
$u^{(j)}_{-1}(0)$.
This will allow us to prove formula
(\ref{subprincipal symbol of OI at time zero}).

Before starting explicit calculations we observe that
equations (\ref{Proof of formula equation 1}) can be equivalently rewritten as
\begin{equation}
\label{Proof of formula equation 3}
P^{(l)}
\,
\left.
\bigl[\mathfrak{S}^{(j)}_{-1}f^{(j)}_1+\mathfrak{S}^{(j)}_0f^{(j)}_0\bigr]
\right|_{t=0}
=0,
\qquad l\ne j,
\end{equation}
where $P^{(l)}:=[v^{(l)}(y,\eta)]\,[v^{(l)}(y,\eta)]^*$
is the orthogonal projection onto the eigenspace corresponding to
the (normalised) eigenvector $v^{(l)}(y,\eta)$ of the principal
symbol.
We will deal with
(\ref{Proof of formula equation 3})
rather than with
(\ref{Proof of formula equation 1}).
This is simply a matter of convenience.

\subsection{Part 1 of the proof of formula (\ref{subprincipal symbol of OI at time zero})}
\label{Part 1}

Our task in this subsection is to calculate
the LHS of (\ref{Proof of formula equation 3}).
In our calculations we use the explicit formula
(\ref{formula for principal symbol of oscillatory integral})
for the principal symbol $u^{(j)}_0(t;y,\eta)$
which was proved in Section~\ref{Leading transport equations}.

At $t=0$ formula (\ref{Leading transport equations equation 4}) reads
\[
\left.
\bigl[\mathfrak{S}^{(j)}_{-1}f^{(j)}_1\bigr]
\right|_{t=0}
=
i
\left.
\left[
\frac{\partial^2}{\partial x^\alpha\eta_\alpha}
\bigl(
A_1(x,\eta)
-h^{(j)}(y,\eta)
-(x-y)^\gamma h^{(j)}_{y^\gamma}(y,\eta)
\bigr)
P^{(j)}(y,\eta)
\right]
\right|_{x=y}
\]
which gives us
\begin{equation}
\label{Proof of formula equation 4}
\left.
\bigl[\mathfrak{S}^{(j)}_{-1}f^{(j)}_1\bigr]
\right|_{t=0}
=
i
\left[
(A_1-h^{(j)})_{y^\alpha\eta_\alpha}P^{(j)}
+
(A_1-h^{(j)})_{y^\alpha}P^{(j)}_{\eta_\alpha}
\right].
\end{equation}
In the latter formula we dropped, for the sake of brevity,
the arguments $(y,\eta)$.

At $t=0$ formula (\ref{Leading transport equations equation 5}) reads
\begin{multline}
\label{Proof of formula equation 5}
\left.
\bigl[\mathfrak{S}^{(j)}_0f^{(j)}_0\bigr]
\right|_{t=0}
=
-i\{v^{(j)},h^{(j)}\}[v^{(j)}]^*
+
\left(
A_0
-
q^{(j)}
+
\frac i2h^{(j)}_{y^\alpha\eta_\alpha}
\right)
P^{(j)}
\\
+[A_1-h^{(j)}]u^{(j)}_{-1}(0)\,,
\end{multline}
where $q^{(j)}$ is the function
(\ref{phase appearing in principal symbol})
and we dropped, for the sake of brevity,
the arguments $(y,\eta)$.
Note that in writing down
(\ref{Proof of formula equation 5})
we used the fact that
\[
\left.
\left[
(d_{\varphi^{(j)}})^{-1}
\partial_t
d_{\varphi^{(j)}}
\right]
\right|_{(t,x;y,\eta)=(0,y;y,\eta)}
=-\frac12
h^{(j)}_{y^\alpha\eta_\alpha}(y,\eta)\,,
\]
compare with formula
(\ref{Leading transport equations equation 20}).

Substituting formulae
(\ref{Proof of formula equation 4})
and
(\ref{Proof of formula equation 5})
into
(\ref{Proof of formula equation 3})
we get
\begin{equation}
\label{Proof of formula equation 6}
(h^{(l)}-h^{(j)})P^{(l)}u^{(j)}_{-1}(0)+P^{(l)}B^{(j)}_0=0,
\qquad l\ne j,
\end{equation}
where
\begin{equation}
\label{Part 1 result}
B^{(j)}_0=
\left(
A_0-q^{(j)}-\frac i2h^{(j)}_{y^\alpha\eta_\alpha}+i(A_1)_{y^\alpha\eta_\alpha}
\right)
P^{(j)}
-i
h^{(j)}_{\eta_\alpha}P^{(j)}_{y^\alpha}
+i(A_1)_{y^\alpha}P^{(j)}_{\eta_\alpha}.
\end{equation}
The subscript in $B^{(j)}_0$ indicates the degree of homogeneity in $\eta$.

\subsection{Part 2 of the proof of formula (\ref{subprincipal symbol of OI at time zero})}
\label{Part 2}

Our task in this subsection is to
solve the linear algebraic system
(\ref{Proof of formula equation 6}),
(\ref{Proof of formula equation 2})
for the unknowns
$u^{(j)}_{-1}(0)$.

It is easy to see that
the unique solution to the system
(\ref{Proof of formula equation 6}),
(\ref{Proof of formula equation 2})
is
\begin{equation}
\label{Part 2 result}
u^{(j)}_{-1}(0)
=\sum_{l\ne j}
\frac
{P^{(l)}B^{(j)}_0+P^{(j)}B^{(l)}_0}
{h^{(j)}-h^{(l)}}\,.
\end{equation}
Summation in (\ref{Part 2 result}) is carried out over all $l$
different from $j$.

\subsection{Part 3 of the proof of formula (\ref{subprincipal symbol of OI at time zero})}
\label{Part 3}

Our task in this subsection is to calculate $[U^{(j)}(0)]_\mathrm{sub}$.

We have
\begin{equation}
\label{subprincipal symbol of Uj0 equation 1}
[U^{(j)}(0)]_\mathrm{sub}
=u^{(j)}_{-1}(0)-\frac i2P^{(j)}_{y^\alpha\eta_\alpha}.
\end{equation}
Here the sign in front of $\frac i2$ is opposite to that in
(\ref{definition of subprincipal symbol})
because the way we write $U^{(j)}(0)$ is using the dual symbol.

Substituting
(\ref{Part 2 result})
and
(\ref{Part 1 result})
into (\ref{subprincipal symbol of Uj0 equation 1})
we get
\begin{multline}
\label{subprincipal symbol of Uj0 equation 2}
[U^{(j)}(0)]_\mathrm{sub}
=
-\frac i2P^{(j)}_{y^\alpha\eta_\alpha}
+\sum_{l\ne j}\frac1{h^{(j)}-h^{(l)}}
\\
\times
\bigl(
P^{(l)}
[
(A_0+i(A_1)_{y^\alpha\eta_\alpha})P^{(j)}
-ih^{(j)}_{\eta_\alpha}P^{(j)}_{y^\alpha}
+i(A_1)_{y^\alpha}P^{(j)}_{\eta_\alpha}
]
\\
\qquad\qquad+
P^{(j)}
[
(A_0+i(A_1)_{y^\alpha\eta_\alpha})P^{(l)}
-ih^{(l)}_{\eta_\alpha}P^{(l)}_{y^\alpha}
+i(A_1)_{y^\alpha}P^{(l)}_{\eta_\alpha}
]
\bigr)
\\
=
\sum_{l\ne j}
\frac
{
P^{(l)}A_\mathrm{sub}P^{(j)}
+
P^{(j)}A_\mathrm{sub}P^{(l)}
}
{
h^{(j)}-h^{(l)}
}
+\frac i2
\Bigl(
-P^{(j)}_{y^\alpha\eta_\alpha}
+
\sum_{l\ne j}
\frac
{
G_{jl}
}
{
h^{(j)}-h^{(l)}
}
\Bigr)\,,
\end{multline}
where
\begin{multline*}
G_{jl}:=
P^{(l)}
[
(A_1)_{y^\alpha\eta_\alpha}P^{(j)}
-2h^{(j)}_{\eta_\alpha}P^{(j)}_{y^\alpha}
+2(A_1)_{y^\alpha}P^{(j)}_{\eta_\alpha}
]
\\
+
P^{(j)}
[
(A_1)_{y^\alpha\eta_\alpha}P^{(l)}
-2h^{(l)}_{\eta_\alpha}P^{(l)}_{y^\alpha}
+2(A_1)_{y^\alpha}P^{(l)}_{\eta_\alpha}
]
\,.
\end{multline*}

We have
\begin{multline*}
G_{jl}
=
2P^{(l)}\{A_1,P^{(j)}\}
+
2P^{(j)}\{A_1,P^{(l)}\}
\\
+
P^{(l)}
[
(A_1-h^{(j)})_{y^\alpha\eta_\alpha}P^{(j)}
+2(A_1-h^{(j)})_{\eta_\alpha}P^{(j)}_{y^\alpha}
]
\\
+
P^{(j)}
[
(A_1-h^{(l)})_{y^\alpha\eta_\alpha}P^{(l)}
+2(A_1-h^{(l)})_{\eta_\alpha}P^{(l)}_{y^\alpha}
]
\\
=
2P^{(l)}\{A_1,P^{(j)}\}
+
2P^{(j)}\{A_1,P^{(l)}\}
-
P^{(l)}\{A_1-h^{(j)},P^{(j)}\}
-
P^{(j)}\{A_1-h^{(l)},P^{(l)}\}
\\
+
P^{(l)}
[
(A_1-h^{(j)})_{y^\alpha\eta_\alpha}P^{(j)}
+(A_1-h^{(j)})_{\eta_\alpha}P^{(j)}_{y^\alpha}
+(A_1-h^{(j)})_{y^\alpha}P^{(j)}_{\eta_\alpha}
]
\\
+
P^{(j)}
[
(A_1-h^{(l)})_{y^\alpha\eta_\alpha}P^{(l)}
+(A_1-h^{(l)})_{\eta_\alpha}P^{(l)}_{y^\alpha}
+(A_1-h^{(l)})_{y^\alpha}P^{(l)}_{\eta_\alpha}
]
\\
=
P^{(l)}\{A_1+h^{(j)},P^{(j)}\}
+
P^{(j)}\{A_1+h^{(l)},P^{(l)}\}
\\
-
P^{(l)}
(A_1-h^{(j)})
P^{(j)}_{y^\alpha\eta_\alpha}
-
P^{(j)}
(A_1-h^{(l)})
P^{(l)}_{y^\alpha\eta_\alpha}
\\
=
P^{(l)}\{A_1+h^{(j)},P^{(j)}\}
+
P^{(j)}\{A_1+h^{(l)},P^{(l)}\}
\\
-
P^{(l)}
(h^{(l)}-h^{(j)})
P^{(j)}_{y^\alpha\eta_\alpha}
-
P^{(j)}
(h^{(j)}-h^{(l)})
P^{(l)}_{y^\alpha\eta_\alpha}
\\
=
P^{(l)}\{A_1+h^{(j)},P^{(j)}\}
+
P^{(j)}\{A_1+h^{(l)},P^{(l)}\}
+(h^{(j)}-h^{(l)})
(
P^{(l)}
P^{(j)}_{y^\alpha\eta_\alpha}
-
P^{(j)}
P^{(l)}_{y^\alpha\eta_\alpha}
)\,,
\end{multline*}
so formula (\ref{subprincipal symbol of Uj0 equation 2}) can be rewritten as
\begin{multline}
\label{subprincipal symbol of Uj0 equation 3}
[U^{(j)}(0)]_\mathrm{sub}
=
\frac i2
\Bigl(
-P^{(j)}_{y^\alpha\eta_\alpha}
+
\sum_{l\ne j}
(
P^{(l)}
P^{(j)}_{y^\alpha\eta_\alpha}
-
P^{(j)}
P^{(l)}_{y^\alpha\eta_\alpha}
)
\Bigr)
\\
+
\frac12
\sum_{l\ne j}
\frac
{
P^{(l)}(2A_\mathrm{sub}P^{(j)}+i\{A_1+h^{(j)},P^{(j)}\})
+
P^{(j)}(2A_\mathrm{sub}P^{(l)}+i\{A_1+h^{(l)},P^{(l)}\})
}
{
h^{(j)}-h^{(l)}
}\,.
\end{multline}

But
\begin{multline*}
\sum_{l\ne j}
(
P^{(l)}
P^{(j)}_{y^\alpha\eta_\alpha}
-
P^{(j)}
P^{(l)}_{y^\alpha\eta_\alpha}
)
=
\Bigl(\,
\sum_{l\ne j}
P^{(l)}
\Bigr)
P^{(j)}_{y^\alpha\eta_\alpha}
-
P^{(j)}
\Bigl(\,
\sum_{l\ne j}
P^{(l)}
\Bigr)_{y^\alpha\eta_\alpha}
\\
=(I-P^{(j)})P^{(j)}_{y^\alpha\eta_\alpha}
-P^{(j)}(I-P^{(j)})_{y^\alpha\eta_\alpha}
=P^{(j)}_{y^\alpha\eta_\alpha},
\end{multline*}
so formula (\ref{subprincipal symbol of Uj0 equation 3}) can be simplified to read
\begin{multline}
\label{Part 3 result}
[U^{(j)}(0)]_\mathrm{sub}
\\
=
\frac12
\sum_{l\ne j}
\frac
{
P^{(l)}(2A_\mathrm{sub}P^{(j)}+i\{A_1+h^{(j)},P^{(j)}\})
+
P^{(j)}(2A_\mathrm{sub}P^{(l)}+i\{A_1+h^{(l)},P^{(l)}\})
}
{
h^{(j)}-h^{(l)}
}
\,.
\end{multline}

\subsection{Part 4 of the proof of formula (\ref{subprincipal symbol of OI at time zero})}
\label{Part 4}

Our task in this subsection is to calculate $\operatorname{tr}[U^{(j)}(0)]_\mathrm{sub}$.

Formula (\ref{Part 3 result}) implies
\begin{equation}
\label{trace of subprincipal symbol of Uj0 equation 1}
\operatorname{tr}[U^{(j)}(0)]_\mathrm{sub}
=
\frac i2\operatorname{tr}\sum_{l\ne j}
\frac
{
P^{(l)}\{A_1,P^{(j)}\}
+
P^{(j)}\{A_1,P^{(l)}\}
}
{
h^{(j)}-h^{(l)}
}
\,.
\end{equation}
Put $A_1=\sum_kh^{(k)}P^{(k)}$ and observe that
\begin{itemize}
\item
terms with the derivatives of $h$ vanish and
\item
the only $k$ which may give nonzero contributions are $k=j$ and $k=l$.
\end{itemize}
Thus, formula (\ref{trace of subprincipal symbol of Uj0 equation 1}) becomes
\begin{multline}
\label{trace of subprincipal symbol of Uj0 equation 2}
\operatorname{tr}[U^{(j)}(0)]_\mathrm{sub}
=
\frac i2\operatorname{tr}\sum_{l\ne j}
\frac1
{
h^{(j)}-h^{(l)}
}
\\
\times\bigl(
h^{(j)}
[
P^{(l)}\{P^{(j)},P^{(j)}\}
+
P^{(j)}\{P^{(j)},P^{(l)}\}
]
+
h^{(l)}
[
P^{(l)}\{P^{(l)},P^{(j)}\}
+
P^{(j)}\{P^{(l)},P^{(l)}\}
]
\bigr).
\end{multline}

We claim that
\begin{multline}
\label{Part 4 auxiliary equation 1}
\operatorname{tr}(P^{(l)}\{P^{(j)},P^{(j)}\})
=
\operatorname{tr}(P^{(j)}\{P^{(j)},P^{(l)}\})
\\
=
-\operatorname{tr}(P^{(l)}\{P^{(l)},P^{(j)}\})
=
-\operatorname{tr}(P^{(j)}\{P^{(l)},P^{(l)}\})
=[v^{(l)}]^*\{v^{(j)},[v^{(j)}]^*\}v^{(l)}
\\
=([v^{(l)}]^*v^{(j)}_{y^\alpha})([v^{(j)}_{\eta_\alpha}]^*v^{(l)})
-([v^{(l)}]^*v^{(j)}_{\eta_\alpha})([v^{(j)}_{y^\alpha}]^*v^{(l)}).
\end{multline}
These facts are established by writing the orthogonal projections
in terms of the eigenvectors and using, if required,
the identities
\[
[v^{(l)}_{y^\alpha}]^*v^{(j)}+[v^{(l)}]^*v^{(j)}_{y^\alpha}=0,
\qquad
[v^{(l)}_{\eta_\alpha}]^*v^{(j)}+[v^{(l)}]^*v^{(j)}_{\eta_\alpha}=0,
\]
\[
[v^{(j)}_{y^\alpha}]^*v^{(l)}+[v^{(j)}]^*v^{(l)}_{y^\alpha}=0,
\qquad
[v^{(j)}_{\eta_\alpha}]^*v^{(l)}+[v^{(j)}]^*v^{(l)}_{\eta_\alpha}=0.
\]
In view of the identities (\ref{Part 4 auxiliary equation 1})
formula (\ref{trace of subprincipal symbol of Uj0 equation 2})
can be rewritten as
\begin{multline}
\label{Part 4 auxiliary equation 2}
\operatorname{tr}[U^{(j)}(0)]_\mathrm{sub}
=
i\operatorname{tr}
\sum_{l\ne j}
P^{(l)}\{P^{(j)},P^{(j)}\}
\\
=
i\operatorname{tr}
(\{P^{(j)},P^{(j)}\}-P^{(j)}\{P^{(j)},P^{(j)}\})
=
-i\operatorname{tr}
(P^{(j)}\{P^{(j)},P^{(j)}\}).
\end{multline}

It remains only to simplify the expression in the RHS of (\ref{Part 4 auxiliary equation 2}).
We have
\begin{multline}
\label{Part 4 auxiliary equation 3}
\operatorname{tr}
(P^{(j)}\{P^{(j)},P^{(j)}\})
=
\{[v^{(j)}]^*,v^{(j)}\}
\\
+[([v^{(j)}]^*v^{(j)}_{y^\alpha})([v^{(j)}]^*v^{(j)}_{\eta_\alpha})-([v^{(j)}]^*v^{(j)}_{\eta_\alpha})([v^{(j)}]^*v^{(j)}_{y^\alpha})]
\\
+[([v^{(j)}_{y^\alpha}]^*v^{(j)})([v^{(j)}_{\eta_\alpha}]^*v^{(j)})-([v^{(j)}_{\eta_\alpha}]^*v^{(j)})([v^{(j)}_{y^\alpha}]^*v^{(j)})]
\\
+[([v^{(j)}]^*v^{(j)}_{y^\alpha})([v^{(j)}_{\eta_\alpha}]^*v^{(j)})-([v^{(j)}]^*v^{(j)}_{\eta_\alpha})([v^{(j)}_{y^\alpha}]^*v^{(j)})]
\\
=
\{[v^{(j)}]^*,v^{(j)}\}
+[([v^{(j)}]^*v^{(j)}_{y^\alpha})([v^{(j)}_{\eta_\alpha}]^*v^{(j)})-([v^{(j)}]^*v^{(j)}_{\eta_\alpha})([v^{(j)}_{y^\alpha}]^*v^{(j)})]
\\
=
\{[v^{(j)}]^*,v^{(j)}\}
-[([v^{(j)}]^*v^{(j)}_{y^\alpha})([v^{(j)}]^*v^{(j)}_{\eta_\alpha})-([v^{(j)}]^*v^{(j)}_{\eta_\alpha})([v^{(j)}]^*v^{(j)}_{y^\alpha})]
\\
=
\{[v^{(j)}]^*,v^{(j)}\}.
\end{multline}
Formulae
(\ref{Part 4 auxiliary equation 2})
and
(\ref{Part 4 auxiliary equation 3})
imply formula (\ref{subprincipal symbol of OI at time zero}).

\section{$\mathrm{U}(1)$ connection}
\label{U(1) connection}

In the preceding Sections
\ref{Algorithm for the construction of the wave group}--\ref{Proof of formula}
we presented technical details
of the construction of the propagator. We saw that
the eigenvectors of the principal symbol, $v^{(j)}(x,\xi)$, play a major role
in this construction. As pointed out in Section~\ref{Main results},
each of these eigenvectors is
defined up to a $\mathrm{U}(1)$ gauge transformation
(\ref{gauge transformation of the eigenvector}),
(\ref{phase appearing in gauge transformation}).
In the end, the full symbols
(\ref{decomposition of symbol of OI into homogeneous components})
of our oscillatory integrals $U^{(j)}(t)$
do not depend on the choice of gauge for the eigenvectors $v^{(j)}(x,\xi)$.
However, the effect of the gauge transformation
(\ref{gauge transformation of the eigenvector}),
(\ref{phase appearing in gauge transformation})
is not as trivial as it may appear at first sight.
We will demonstrate in this section that the gauge transformation
(\ref{gauge transformation of the eigenvector}),
(\ref{phase appearing in gauge transformation})
shows up, in the form of invariantly defined curvature, in the lower
order terms $u^{(j)}_{-1}(t;y,\eta)$ of the symbols of our oscillatory integrals $U^{(j)}(t)$.
More precisely, we will show that the RHS of
formula~(\ref{subprincipal symbol of OI at time zero})
is the scalar curvature of a connection associated with the gauge transformation
(\ref{gauge transformation of the eigenvector}),
(\ref{phase appearing in gauge transformation}).
Further on in this section, until the very last paragraph, the index $j$ enumerating eigenvalues and
eigenvectors of the principal symbol is assumed to be fixed.

Consider a smooth curve $\Gamma\subset T'M$ connecting points $(y,\eta)$ and $(x,\xi)$.
We write this curve in parametric form as $(z(t),\zeta(t))$, $t\in[0,1]$,
so that $(z(0),\zeta(0))=(y,\eta)$ and $(z(1),\zeta(1))=(x,\xi)$.
Put
\begin{equation}
\label{derivative of eigenvector is orthogonal to eigenvector auxiliary}
w(t):=e^{i\phi(t)}v^{(j)}(z(t),\zeta(t))\,,
\end{equation}
where
$\phi:[0,1]\to\mathbb{R}$
is an unknown function which is to be determined from the condition
\begin{equation}
\label{derivative of eigenvector is orthogonal to eigenvector}
iw^*\dot w=0
\end{equation}
with the dot indicating the derivative with respect to the parameter $t$.
Substituting (\ref{derivative of eigenvector is orthogonal to eigenvector auxiliary})
into
(\ref{derivative of eigenvector is orthogonal to eigenvector})
we get an ordinary differential equation for $\phi$ which
is easily solved, giving
\begin{multline}
\label{formula for phi(1)}
\phi(1)
=\phi(0)+\int_0^1(\dot z^\alpha(t)\,P_\alpha(z(t),\zeta(t))+\dot\zeta_\gamma(t)\,Q^\gamma(z(t),\zeta(t)))\,dt
\\
=\phi(0)+\int_\Gamma(P_\alpha dz^\alpha+Q^\gamma d\zeta_\gamma)\,,
\end{multline}
where
\begin{equation}
\label{formula for P and Q}
P_\alpha:=i[v^{(j)}]^*v^{(j)}_{z^\alpha},
\qquad
Q^\gamma:=i[v^{(j)}]^*v^{(j)}_{\zeta_\gamma}.
\end{equation}
Note that the $2n$-component real quantity $(P_\alpha,Q^\gamma)$
is a covector field (1-form) on $T'M$. This quantity already appeared
in Section~\ref{Main results} as formula (\ref{electromagnetic covector potential}).

Put $f(y,\eta):=e^{i\phi(0)}$, $f(x,\xi):=e^{i\phi(1)}$
and rewrite formula (\ref{formula for phi(1)}) as
\begin{equation}
\label{formula for a(1)}
f(x,\xi)
=f(y,\eta)\,e^{i\int_\Gamma(P_\alpha dz^\alpha+Q^\gamma d\zeta_\gamma)}.
\end{equation}
Let us identify the group $\mathrm{U}(1)$ with the unit circle in the complex
plane, i.e. with $f\in\mathbb{C}$, $|f|=1$.
We see that formulae (\ref{formula for a(1)}) and (\ref{formula for P and Q})
give us a rule for the parallel transport of elements
of the group $\mathrm{U}(1)$ along curves in $T'M$. This is the natural
$\mathrm{U}(1)$ connection generated by the normalised field of columns of
complex-valued scalars
\begin{equation}
\label{jth eigenvector of the principal symbol}
v^{(j)}(z,\zeta)=
\bigl(
\begin{matrix}v^{(j)}_1(z,\zeta)&\ldots&v^{(j)}_m(z,\zeta)\end{matrix}
\bigr)^T.
\end{equation}
Recall that the $\Gamma$ appearing in formula (\ref{formula for a(1)}) is a curve
connecting points $(y,\eta)$ and $(x,\xi)$, whereas
the $v^{(j)}(z,\zeta)$ appearing in formulae
(\ref{formula for P and Q}) and (\ref{jth eigenvector of the principal symbol})
enters our construction as
an eigenvector of the principal symbol of our $m\times m$ matrix pseudo\-differential
operator $A$.

In practice, dealing with a connection is not as convenient as dealing with
the covariant derivative $\nabla$. The covariant derivative
corresponding to the connection (\ref{formula for a(1)}) is determined as follows.
Let us view the $(x,\xi)$ appearing in formula (\ref{formula for a(1)})
as a variable which takes values close to $(y,\eta)$,
and suppose that the curve $\Gamma$ is a short straight (in local coordinates)
line segment connecting the point $(y,\eta)$ with the point $(x,\xi)$.
We want the covariant derivative of our function
$f(x,\xi)$, evaluated at $(y,\eta)$, to be zero.
Examination of formula (\ref{formula for a(1)}) shows that
the unique covariant derivative satisfying this condition is
\begin{equation}
\label{formula for U(1) covariant derivative}
\nabla_\alpha:=\partial/\partial x^\alpha-iP_\alpha(x,\xi),
\qquad
\nabla^\gamma:=\partial/\partial\xi_\gamma-iQ^\gamma(x,\xi).
\end{equation}

We define the curvature of our $\mathrm{U}(1)$ connection as
\begin{equation}
\label{definition of U(1) curvature}
R:=
-i
\begin{pmatrix}
\nabla_\alpha\nabla_\beta-\nabla_\beta\nabla_\alpha&
\nabla_\alpha\nabla^\delta-\nabla^\delta\nabla_\alpha
\\
\nabla^\gamma\nabla_\beta-\nabla_\beta\nabla^\gamma&
\nabla^\gamma\nabla^\delta-\nabla^\delta\nabla^\gamma
\end{pmatrix}.
\end{equation}
It may seem that the entries of the $(2n)\times(2n)$ matrix (\ref{definition of U(1) curvature})
are differential operators. They are, in fact, operators of multiplication
by ``scalar functions''. Namely, the more explicit form of (\ref{definition of U(1) curvature}) is
\begin{equation}
\label{explicit formula for U(1) curvature}
R=
\begin{pmatrix}
\frac{\partial P_\alpha}{\partial x^\beta}-\frac{\partial P_\beta}{\partial x^\alpha}&
\frac{\partial P_\alpha}{\partial\xi_\delta}-\frac{\partial Q^\delta}{\partial x^\alpha}
\\
\frac{\partial Q^\gamma}{\partial x^\beta}-\frac{\partial P_\beta}{\partial\xi_\gamma}&
\frac{\partial Q^\gamma}{\partial\xi_\delta}-\frac{\partial Q^\delta}{\partial\xi_\gamma}
\end{pmatrix}.
\end{equation}
The $(2n)\times(2n)$\,-\,component real quantity (\ref{explicit formula for U(1) curvature})
is a rank 2 covariant antisymmetric tensor (2-form) on $T'M$.
It is an analogue of the electromagnetic tensor.

Substituting (\ref{formula for P and Q}) into
(\ref{explicit formula for U(1) curvature})
we get an expression for curvature in terms of the eigenvector
of the principal symbol
\begin{equation}
\label{more explicit formula for U(1) curvature}
R=i
\begin{pmatrix}
[v^{(j)}_{x^\beta}]^*v^{(j)}_{x^\alpha}-[v^{(j)}_{x^\alpha}]^*v^{(j)}_{x^\beta}&
[v^{(j)}_{\xi_\delta}]^*v^{(j)}_{x^\alpha}-[v^{(j)}_{x^\alpha}]^*v^{(j)}_{\xi_\delta}
\\
[v^{(j)}_{x^\beta}]^*v^{(j)}_{\xi_\gamma}-[v^{(j)}_{\xi_\gamma}]^*v^{(j)}_{x^\beta}&
[v^{(j)}_{\xi_\delta}]^*v^{(j)}_{\xi_\gamma}-[v^{(j)}_{\xi_\gamma}]^*v^{(j)}_{\xi_\delta}
\end{pmatrix}.
\end{equation}
Examination of formula (\ref{more explicit formula for U(1) curvature}) shows that,
as expected, curvature is invariant under the gauge transformation
(\ref{gauge transformation of the eigenvector}),
(\ref{phase appearing in gauge transformation}).

It is natural to take the trace of the upper right block
in (\ref{definition of U(1) curvature}) which,
in the notation (\ref{Poisson bracket on matrix-functions}), gives us
\begin{equation}
\label{scalar curvature of U(1) connection}
-i(\nabla_\alpha\nabla^\alpha-\nabla^\alpha\nabla_\alpha)
=-i\{[v^{(j)}]^*,v^{(j)}\}.
\end{equation}
Thus, we have shown that the RHS of
formula~(\ref{subprincipal symbol of OI at time zero})
is the scalar curvature of our $\mathrm{U}(1)$ connection.

\

We end this section by proving, as promised in Section~\ref{Main results},
formula (\ref{sum of curvatures is zero}) without referring to microlocal analysis.
In the following arguments we use our standard notation for the orthogonal
projections onto the eigenspaces of the principal symbol,
i.e.~we write $P^{(k)}:=v^{(k)}[v^{(k)}]^*$.
We have $\operatorname{tr}\{P^{(j)},P^{(j)}\}=0$
and $\sum_lP^{(l)}=I$
which implies
\begin{multline}
\label{sum of curvatures is zero proof equation 1}
0=\sum_{l,j}\operatorname{tr}(P^{(l)}\{P^{(j)},P^{(j)}\})
\\
=\sum_j\operatorname{tr}(P^{(j)}\{P^{(j)},P^{(j)}\})
+\sum_{l,j:\ l\ne j}\operatorname{tr}(P^{(l)}\{P^{(j)},P^{(j)}\}).
\end{multline}
But, according to formula (\ref{Part 4 auxiliary equation 1}),
for $l\ne j$ we have
\[
\operatorname{tr}(P^{(l)}\{P^{(j)},P^{(j)}\})
=-\operatorname{tr}(P^{(j)}\{P^{(l)},P^{(l)}\}),
\]
so the expression in the last sum in the RHS of (\ref{sum of curvatures is zero proof equation 1})
is antisymmetric in the indices $l,j$, which implies that this sum is zero.
Hence, formula
(\ref{sum of curvatures is zero proof equation 1}) can be rewritten as
$\sum\limits_j\operatorname{tr}(P^{(j)}\{P^{(j)},P^{(j)}\})=0$.
It remains only to note that,
according to formula (\ref{Part 4 auxiliary equation 3}),
$\operatorname{tr}(P^{(j)}\{P^{(j)},P^{(j)}\})=\{[v^{(j)}]^*,v^{(j)}\}$.

\section{Singularity of the propagator at $t=0$}
\label{Singularity of the wave group at time zero}

Following the notation of \cite{mybook}, we denote by
\[
\mathcal{F}_{\lambda\to t}[f(\lambda)]=\hat f(t)=\int e^{-it\lambda}f(\lambda)\,d\lambda
\]
the one-dimensional Fourier transform and by
\[
\mathcal{F}^{-1}_{t\to\lambda}[\hat f(t)]=f(\lambda)=(2\pi)^{-1}\int e^{it\lambda}\hat f(t)\,dt
\]
its inverse.

Suppose that we have a Hamiltonian trajectory
$(x^{(j)}(t;y,\eta),\xi^{(j)}(t;y,\eta))$
and a real number $T>0$ such that
$x^{(j)}(T;y,\eta)=y$. We will say in this case
that we have a loop of length $T$ originating
from the point $y\in M$.

\begin{remark}
\label{remark on reversibility}
There is no need to consider loops of negative length $T$ because,
given a $T>0$, we have
$x^{(j)}(T;y,\eta^+)=y$
for some $\eta^+\in T'_yM$ if and only if we have
$x^{(j)}(-T;y,\eta^-)=y$
for some $\eta^-\in T'_yM$. Indeed,
it suffices to relate the
$\eta^\pm$ in accordance with
$\eta^\mp=\xi^{(j)}(\pm T;y,\eta^\pm)$.
\end{remark}

Denote by $\mathcal{T}^{(j)}\subset\mathbb{R}$ the set of lengths $T>0$
of all possible loops generated by the Hamiltonian $h^{(j)}$.
Here ``all possible'' refers to all possible starting points
$(y,\eta)\in T'M$ of Hamiltonian trajectories.
It is easy to see that $0\not\in\overline{\mathcal{T}^{(j)}}$.
We put
\[
\mathbf{T}^{(j)}:=
\begin{cases}
\inf\mathcal{T}^{(j)}\quad&\text{if}\quad\mathcal{T}^{(j)}\ne\emptyset,
\\
+\infty\quad&\text{if}\quad\mathcal{T}^{(j)}=\emptyset.
\end{cases}
\]

In the Riemannian case (i.e.~the case when the Hamiltonian
is a square root of a quadratic polynomial in $\xi$) it is known \cite{sabourau,rotman}
that there is a loop originating from every point of the
manifold $M$ and, moreover, there is an explicit estimate from above for
the number $\mathbf{T}^{(j)}$.
We are not aware of similar results for general Hamiltonians.

We also define
$\mathbf{T}:=\min\limits_{j=1,\ldots,m^+}\mathbf{T}^{(j)}$.

\begin{remark}
\label{remark on negative Hamiltonians}
Note that negative eigenvalues of the principal symbol,
i.e.~Hamiltonians $h^{(j)}(x,\xi)$ with negative index
$j=-1,\ldots,-m^-$,
do not affect the asymptotic formulae
we are about to derive. This is because we are dealing
with the case $\lambda\to+\infty$ rather than $\lambda\to-\infty$.
\end{remark}

Denote by
\begin{equation}
\label{definition of integral kernel of wave group}
u(t,x,y):=
\sum_k e^{-it\lambda_k}v_k(x)[v_k(y)]^*
\end{equation}
the integral kernel of the propagator (\ref{definition of wave group}).
The quantity (\ref{definition of integral kernel of wave group})
can be understood as a distribution in the variable
$t\in\mathbb{R}$ depending on the parameters $x,y\in M$.

The main result of this section is the following
\begin{lemma}
\label{Singularity of the wave group at time zero lemma}
Let $\hat\rho:\mathbb{R}\to\mathbb{C}$ be an infinitely smooth function such that
\begin{equation}
\label{condition on hat rho 1}
\operatorname{supp}\hat\rho\subset(-\mathbf{T},\mathbf{T}),
\end{equation}
\begin{equation}
\label{condition on hat rho 2}
\hat\rho(0)=1,
\end{equation}
\begin{equation}
\label{condition on hat rho 3}
\hat\rho'(0)=0.
\end{equation}
Then, uniformly over $y\in M$, we have
\begin{equation}
\label{Singularity of the wave group at time zero lemma formula}
\mathcal{F}^{-1}_{t\to\lambda}[\hat\rho(t)\operatorname{tr}u(t,y,y)]=
n\,a(y)\,\lambda^{n-1}+(n-1)\,b(y)\,\lambda^{n-2}+O(\lambda^{n-3})
\end{equation}
as $\lambda\to+\infty$.
The densities $a(y)$ and $b(y)$ appearing in the RHS of formula
(\ref{Singularity of the wave group at time zero lemma formula})
are defined in accordance with formulae
(\ref{formula for a(x)}) and (\ref{formula for b(x)}).
\end{lemma}

\emph{Proof\ }
Denote by $(S^*_yM)^{(j)}$ the $(n-1)$-dimensional unit cosphere in the cotangent fibre
defined by the equation $h^{(j)}(y,\eta)=1$
and denote by $d(S^*_yM)^{(j)}$ the
surface area element on $(S^*_yM)^{(j)}$
defined by the condition
$d\eta=d(S^*_yM)^{(j)}\,dh^{(j)}$.
The latter means that we use spherical coordinates in the cotangent fibre
with the Hamiltonian $h^{(j)}$
playing the role of the radial coordinate, see subsection 1.1.10 of \cite{mybook} for details.
In particular, as explained in subsection 1.1.10 of \cite{mybook},
our surface area element $d(S^*_yM)^{(j)}$ is expressed via the Euclidean surface area element as
\[
d(S^*_yM)^{(j)}=
\biggl(\,\sum_{\alpha=1}^n\bigl(h^{(j)}_{\eta_\alpha}(y,\eta)\bigr)^2\biggr)^{-1/2}
\times\,
\text{Euclidean surface area element}
\,.
\]
Denote also
$\,{d{\hskip-1pt\bar{}}\hskip1pt}(S^*_yM)^{(j)}:=(2\pi)^{-n}\,d(S^*_yM)^{(j)}\,$.

According to Corollary 4.1.5 from \cite{mybook} we have
uniformly over $y\in M$
\begin{multline}
\label{Singularity of the wave group at time zero lemma equation 1}
\mathcal{F}^{-1}_{t\to\lambda}[\hat\rho(t)\operatorname{tr}u(t,y,y)]=
\\
\sum_{j=1}^{m^+}
\left(c^{(j)}(y)\,\lambda^{n-1}+d^{(j)}(y)\,\lambda^{n-2}+e^{(j)}(y)\,\lambda^{n-2}\right)
+O(\lambda^{n-3})\,,
\end{multline}
where
\begin{equation}
\label{Singularity of the wave group at time zero lemma equation 2}
c^{(j)}(y)=\int\limits_{(S^*_yM)^{(j)}}
\operatorname{tr}u^{(j)}_0(0;y,\eta)
\,{d{\hskip-1pt\bar{}}\hskip1pt}(S^*_yM)^{(j)}\,,
\end{equation}
\begin{multline}
\label{Singularity of the wave group at time zero lemma equation 3}
d^{(j)}(y)=
\\
(n-1)\int\limits_{(S^*_yM)^{(j)}}
\operatorname{tr}
\left(
-\,i\,\dot u^{(j)}_0(0;y,\eta)
+\frac i2\bigl\{u^{(j)}_0\bigr|_{t=0}\,,h^{(j)}\bigr\}(y,\eta)
\right)
{d{\hskip-1pt\bar{}}\hskip1pt}(S^*_yM)^{(j)}\,,
\end{multline}
\begin{equation}
\label{Singularity of the wave group at time zero lemma equation 4}
e^{(j)}(y)=\int\limits_{(S^*_yM)^{(j)}}
\operatorname{tr}[U^{(j)}(0)]_\mathrm{sub}(y,\eta)
\,{d{\hskip-1pt\bar{}}\hskip1pt}(S^*_yM)^{(j)}\,.
\end{equation}
Here $u^{(j)}_0(t;y,\eta)$ is the principal symbol of the oscillatory integral
(\ref{algorithm equation 1}) and $\dot u^{(j)}_0(t;y,\eta)$ is its time derivative.
Note that in writing the term with the Poisson bracket in
(\ref{Singularity of the wave group at time zero lemma equation 3})
we took account of the fact that Poisson brackets in \cite{mybook}
and in the current paper have opposite signs.

Observe that the integrands in formulae
(\ref{Singularity of the wave group at time zero lemma equation 2})
and
(\ref{Singularity of the wave group at time zero lemma equation 3})
are positively homogeneous in $\eta$ of degree 0,
whereas the integrand in formula
(\ref{Singularity of the wave group at time zero lemma equation 4})
is positively homogeneous in $\eta$ of degree $-1$.
In order to have the same degree of homogeneity,  we rewrite
formula
(\ref{Singularity of the wave group at time zero lemma equation 4})
in equivalent form
\begin{equation}
\label{Singularity of the wave group at time zero lemma equation 5}
e^{(j)}(y)=\int\limits_{(S^*_yM)^{(j)}}
\bigl(
h^{(j)}\operatorname{tr}[U^{(j)}(0)]_\mathrm{sub}
\bigr)
(y,\eta)
\,{d{\hskip-1pt\bar{}}\hskip1pt}(S^*_yM)^{(j)}\,.
\end{equation}

Switching from surface integrals to volume integrals with the help of formula (1.1.15) from \cite{mybook},
we rewrite formulae
(\ref{Singularity of the wave group at time zero lemma equation 2}),
(\ref{Singularity of the wave group at time zero lemma equation 3})
and
(\ref{Singularity of the wave group at time zero lemma equation 5})
as
\begin{equation}
\label{Singularity of the wave group at time zero lemma equation 6}
c^{(j)}(y)=n\int\limits_{h^{(j)}(y,\eta)<1}
\operatorname{tr}u^{(j)}_0(0;y,\eta)
\,{d{\hskip-1pt\bar{}}\hskip1pt}\eta\,,
\end{equation}
\begin{multline}
\label{Singularity of the wave group at time zero lemma equation 7}
d^{(j)}(y)=n(n-1)\times
\\
\int\limits_{h^{(j)}(y,\eta)<1}
\operatorname{tr}
\left(
-\,i\,\dot u^{(j)}_0(0;y,\eta)
+\frac i2\bigl\{u^{(j)}_0\bigr|_{t=0}\,,h^{(j)}\bigr\}(y,\eta)
\right)
{d{\hskip-1pt\bar{}}\hskip1pt}\eta\,,
\end{multline}
\begin{equation}
\label{Singularity of the wave group at time zero lemma equation 8}
e^{(j)}(y)=n\int\limits_{h^{(j)}(y,\eta)<1}
\bigl(
h^{(j)}\operatorname{tr}[U^{(j)}(0)]_\mathrm{sub}
\bigr)
(y,\eta)
\,{d{\hskip-1pt\bar{}}\hskip1pt}\eta\,.
\end{equation}

Substituting formulae
(\ref{formula for principal symbol of oscillatory integral})
and
(\ref{phase appearing in principal symbol})
into formulae
(\ref{Singularity of the wave group at time zero lemma equation 6})
and
(\ref{Singularity of the wave group at time zero lemma equation 7})
we get
\begin{equation}
\label{Singularity of the wave group at time zero lemma equation 9}
c^{(j)}(y)=n\int\limits_{h^{(j)}(y,\eta)<1}
{d{\hskip-1pt\bar{}}\hskip1pt}\eta\,,
\end{equation}
\begin{multline}
\label{Singularity of the wave group at time zero lemma equation 10}
d^{(j)}(y)=-n(n-1)\times
\\
\int\limits_{h^{(j)}(y,\eta)<1}
\left(
[v^{(j)}]^*A_\mathrm{sub}v^{(j)}
-\frac i2
\{
[v^{(j)}]^*,A_1-h^{(j)},v^{(j)}
\}
\right)(y,\eta)
\,{d{\hskip-1pt\bar{}}\hskip1pt}\eta\,.
\end{multline}
Substituting formula
(\ref{subprincipal symbol of OI at time zero})
into formula
(\ref{Singularity of the wave group at time zero lemma equation 8})
we get
\begin{equation}
\label{Singularity of the wave group at time zero lemma equation 11}
e^{(j)}(y)=-n\,i\int\limits_{h^{(j)}(y,\eta)<1}
\bigl(
h^{(j)}\{[v^{(j)}]^*,v^{(j)}\}
\bigr)
(y,\eta)
\,{d{\hskip-1pt\bar{}}\hskip1pt}\eta\,.
\end{equation}

Substituting formulae
(\ref{Singularity of the wave group at time zero lemma equation 9})--(\ref{Singularity of the wave group at time zero lemma equation 11})
into formula
(\ref{Singularity of the wave group at time zero lemma equation 1})
we arrive
at (\ref{Singularity of the wave group at time zero lemma formula}).~$\square$

\begin{remark}
The proof of Lemma~\ref{Singularity of the wave group at time zero lemma}
given above was based on the use of Corollary 4.1.5 from \cite{mybook}.
In the actual statement of Corollary 4.1.5 in \cite{mybook}
uniformity in $y\in M$ was not mentioned because the authors were dealing with
a manifold with a boundary. Uniformity reappeared in the subsequent
Theorem 4.2.1 which involved pseudodifferential cut-offs
separating the point $\,y\,$ from the boundary.
\end{remark}

\section{Mollified spectral asymptotics}
\label{Mollified spectral asymptotics}

\begin{theorem}
\label{theorem spectral function mollified}
Let $\rho:\mathbb{R}\to\mathbb{C}$ be a function from Schwartz space $\mathcal{S}(\mathbb{R})$
whose Fourier transform $\hat\rho$ satisfies conditions
(\ref{condition on hat rho 1})--(\ref{condition on hat rho 3}).
Then, uniformly over $x\in M$, we have
\begin{equation}
\label{theorem spectral function mollified formula}
\int e(\lambda-\mu,x,x)\,\rho(\mu)\,d\mu=
a(x)\,\lambda^n+b(x)\,\lambda^{n-1}+
\begin{cases}
O(\lambda^{n-2})\quad&\text{if}\quad{n\ge3},
\\
O(\ln\lambda)\quad&\text{if}\quad{n=2},
\end{cases}
\end{equation}
as $\lambda\to+\infty$.
The densities $a(x)$ and $b(x)$ appearing in the RHS of formula
(\ref{theorem spectral function mollified formula})
are defined in accordance with formulae
(\ref{formula for a(x)}) and (\ref{formula for b(x)}).
\end{theorem}

\emph{Proof\ }
Our spectral function $e(\lambda,x,x)$ was initially defined only for $\lambda>0$,
see formula (\ref{definition of spectral function}). We extend the definition
to the whole real line by setting
\[
e(\lambda,x,x):=0\quad\text{for}\quad\lambda\le0.
\]

Denote by $e'(\lambda,x,x)$ the derivative, with respect to the spectral
parameter, of the spectral function. Here ``derivative'' is understood in the
sense of distributions. The explicit formula for $e'(\lambda,x,x)$ is
\begin{equation}
\label{theorem spectral function mollified equation 2}
e'(\lambda,x,x):=\sum_{k=1}^{+\infty}\|v_k(x)\|^2\,\delta(\lambda-\lambda_k).
\end{equation}

Formula (\ref{theorem spectral function mollified equation 2}) gives us
\begin{equation}
\label{theorem spectral function mollified equation 3}
\int e'(\lambda-\mu,x,x)\,\rho(\mu)\,d\mu=
\sum_{k=1}^{+\infty}\|v_k(x)\|^2\,\rho(\lambda-\lambda_k).
\end{equation}
Formula (\ref{theorem spectral function mollified equation 3})
implies, in particular, that, uniformly over $x\in M$, we have
\begin{equation}
\label{theorem spectral function mollified equation 4}
\int e'(\lambda-\mu,x,x)\,\rho(\mu)\,d\mu=O(|\lambda|^{-\infty})
\quad\text{as}\quad\lambda\to-\infty\,,
\end{equation}
where $O(|\lambda|^{-\infty})$ is shorthand for ``tends to zero faster
than any given inverse power of $|\lambda|$''.

Formula (\ref{theorem spectral function mollified equation 3})
can also be rewritten as
\begin{equation}
\label{theorem spectral function mollified equation 5}
\int e'(\lambda-\mu,x,x)\,\rho(\mu)\,d\mu=
\mathcal{F}^{-1}_{t\to\lambda}[\hat\rho(t)\operatorname{tr}u(t,x,x)]
-\sum_{k\le0}\|v_k(x)\|^2\,\rho(\lambda-\lambda_k)\,,
\end{equation}
where the distribution $u(t,x,y)$ is defined in accordance with
formula (\ref{definition of integral kernel of wave group}).
Clearly, we have
\begin{equation}
\label{theorem spectral function mollified equation 6}
\sum_{k\le0}\|v_k(x)\|^2\,\rho(\lambda-\lambda_k)=O(\lambda^{-\infty})
\quad\text{as}\quad\lambda\to+\infty\,.
\end{equation}
Formulae
(\ref{theorem spectral function mollified equation 5}),
(\ref{theorem spectral function mollified equation 6})
and Lemma~\ref{Singularity of the wave group at time zero lemma}
imply that, uniformly over $x\in M$, we have
\begin{multline}
\label{theorem spectral function mollified equation 7}
\int e'(\lambda-\mu,x,x)\,\rho(\mu)\,d\mu=
\\
n\,a(x)\,\lambda^{n-1}+(n-1)\,b(x)\,\lambda^{n-2}+O(\lambda^{n-3})
\quad\text{as}\quad\lambda\to+\infty\,.
\end{multline}

It remains to note that
\begin{equation}
\label{theorem spectral function mollified equation 8}
\frac d{d\lambda}\int e(\lambda-\mu,x,x)\,\rho(\mu)\,d\mu
=
\int e'(\lambda-\mu,x,x)\,\rho(\mu)\,d\mu\,.
\end{equation}
Formulae
(\ref{theorem spectral function mollified equation 8}),
(\ref{theorem spectral function mollified equation 4})
and
(\ref{theorem spectral function mollified equation 7})
imply
(\ref{theorem spectral function mollified formula}).~$\square$

\begin{theorem}
\label{theorem counting function mollified}
Let $\rho:\mathbb{R}\to\mathbb{C}$ be a function from Schwartz space $\mathcal{S}(\mathbb{R})$
whose Fourier transform $\hat\rho$ satisfies conditions
(\ref{condition on hat rho 1})--(\ref{condition on hat rho 3}).
Then we have
\begin{equation}
\label{theorem counting function mollified formula}
\int N(\lambda-\mu)\,\rho(\mu)\,d\mu=
a\,\lambda^n+b\,\lambda^{n-1}+
\begin{cases}
O(\lambda^{n-2})\quad&\text{if}\quad{n\ge3},
\\
O(\ln\lambda)\quad&\text{if}\quad{n=2},
\end{cases}
\end{equation}
as $\lambda\to+\infty$.
The constants $a$ and $b$ appearing in the RHS of formula
(\ref{theorem counting function mollified formula})
are defined in accordance with formulae
(\ref{a via a(x)}),
(\ref{formula for a(x)}),
(\ref{b via b(x)})
and
(\ref{formula for b(x)}).
\end{theorem}

\emph{Proof\ }
Formula
(\ref{theorem counting function mollified formula})
follows from formula
(\ref{theorem spectral function mollified formula})
by integration over $M$,
see also formula (\ref{definition of counting function}).~$\square$

\

In stating Theorems
\ref{theorem spectral function mollified}
and
\ref{theorem counting function mollified}
we assumed the mollifier $\rho$ to be complex-valued.
This was done for the sake of generality but may seem
unnatural when mollifying real-valued functions
$e(\lambda,x,x)$ and $N(\lambda)$. One can make our
construction look more natural by dealing only with
real-valued mollifiers $\rho$. Note that if the function $\rho$
is real-valued and even then its Fourier transform
$\hat\rho$ is also real-valued and even and, moreover,
condition (\ref{condition on hat rho 3}) is automatically satisfied.

\section{Unmollified spectral asymptotics}
\label{Unmollified spectral asymptotics}

In this section we derive asymptotic formulae for
the spectral function $e(\lambda,x,x)$ and the
counting function $N(\lambda)$ without mollification.
The section is split into two subsections: in the first
we derive one-term asymptotic formulae and
in the second --- two-term asymptotic formulae.

\subsection{One-term spectral asymptotics}
\label{One-term spectral asymptotics}

\begin{theorem}
\label{theorem spectral function unmollified one term}
We have, uniformly over $x\in M$,
\begin{equation}
\label{theorem spectral function unmollified one term formula}
e(\lambda,x,x)=a(x)\,\lambda^n+O(\lambda^{n-1})
\end{equation}
as $\lambda\to+\infty$.
\end{theorem}

\emph{Proof\ }
The result in question is an immediate consequence of
formulae
(\ref{theorem spectral function mollified equation 8}),
(\ref{theorem spectral function mollified equation 7})
and
Theorem~\ref{theorem spectral function mollified}
from the current paper
and Corollary~B.2.2 from \cite{mybook}.~$\square$

\begin{theorem}
\label{theorem counting function unmollified one term}
We have
\begin{equation}
\label{theorem counting function unmollified one term formula}
N(\lambda)=a\lambda^n+O(\lambda^{n-1})
\end{equation}
as $\lambda\to+\infty$.
\end{theorem}

\emph{Proof\ }
Formula
(\ref{theorem counting function unmollified one term formula})
follows from formula
(\ref{theorem spectral function unmollified one term formula})
by integration over $M$,
see also formula (\ref{definition of counting function}).~$\square$

\subsection{Two-term spectral asymptotics}
\label{Two-term spectral asymptotics}

Up till now, in Section~\ref{Mollified spectral asymptotics}
and subsection~\ref{One-term spectral asymptotics},
our logic was to derive asymptotic formulae for the spectral
function $e(\lambda,x,x)$ first and then obtain corresponding
asymptotic formulae for the counting function $N(\lambda)$
by integration over $M$. Such an approach will not work
for two-term asymptotics because
the geometric conditions required for the existence of
two-term asymptotics of $e(\lambda,x,x)$ and $N(\lambda)$
will be different:
for $e(\lambda,x,x)$ the appropriate geometric conditions
will be formulated in terms of \emph{loops},
whereas
for $N(\lambda)$ the appropriate geometric conditions
will be formulated in terms of \emph{periodic trajectories}.

Hence, in this subsection we deal with
the spectral function $e(\lambda,x,x)$
and the counting function $N(\lambda)$ separately.

In what follows the point $y\in M$ is assumed to be fixed.

Denote by $\Pi_y^{(j)}$ the set of normalised ($h^{(j)}(y,\eta)=1$)
covectors $\eta$ which serve as starting points for loops generated by the Hamiltonian
$h^{(j)}$. Here ``starting point'' refers to the starting point
of a Hamiltonian trajectory
$(x^{(j)}(t;y,\eta),\xi^{(j)}(t;y,\eta))$
moving forward in time ($t>0$),
see also Remark~\ref{remark on reversibility}.

The reason we are not interested in large negative $t$ is that the refined
Fourier Tauberian theorem we will be applying,
Theorem~B.5.1 from \cite{mybook},
does not require information regarding large negative $t$.
And the underlying reason for the latter is the fact that the function
we are studying, $e(\lambda,x,x)$ (and, later, $N(\lambda)$), is real-valued.
The real-valuedness of the function $e(\lambda,x,x)$ implies that its
Fourier transform, $\hat e(t,x,x)$, possesses the symmetry
$\hat e(-t,x,x)=\overline{\hat e(t,x,x)}$.

The set $\Pi_y^{(j)}$ is a subset of the $(n-1)$-dimensional unit cosphere
$(S^*_yM)^{(j)}$ and the latter is equipped with
a natural Lebesgue measure, see proof of
Lemma~\ref{Singularity of the wave group at time zero lemma}.
It is known, see Lemma 1.8.2 in \cite{mybook}, that the
set $\Pi_y^{(j)}$ is measurable.

\begin{definition}
\label{definition of nonfocal point 1}
A point $y\in M$ is said to be \emph{nonfocal} if for each
$j=1,\ldots,m^+$ the set $\Pi_y^{(j)}$ has measure zero.
\end{definition}

With regards to the range of the index $j$ in
Definition~\ref{definition of nonfocal point 1},
as well as in sub\-sequent
Definitions~\ref{definition of nonfocal point 2}--\ref{definition of nonperiodicity condition 2},
see Remark~\ref{remark on negative Hamiltonians}.

We call a loop of length $T>0$ \emph{absolutely focused} if
the function
\[
|x^{(j)}(T;y,\eta)-y|^2
\]
has an infinite order zero in the variable $\eta$, and we denote
by $(\Pi_y^a)^{(j)}$ the set of normalised ($h^{(j)}(y,\eta)=1$)
covectors $\eta$ which serve as starting points for absolutely focused loops
generated by the Hamiltonian $h^{(j)}$.
It is known, see Lemma 1.8.3 in \cite{mybook}, that the
set $(\Pi_y^a)^{(j)}$ is measurable and,
moreover, the set $\Pi_y^{(j)}\setminus(\Pi_y^a)^{(j)}$ has measure zero.
This allows us to reformulate Definition~\ref{definition of nonfocal point 1}
as follows.

\begin{definition}
\label{definition of nonfocal point 2}
A point $y\in M$ is said to be \emph{nonfocal} if for each
$j=1,\ldots,m^+$ the set $(\Pi_y^a)^{(j)}$ has measure zero.
\end{definition}

In practical applications it is easier to work with
Definition~\ref{definition of nonfocal point 2}
because the set $(\Pi_y^a)^{(j)}$ is usually much
thinner than the set $\Pi_y^{(j)}$.

In order to derive a two-term asymptotic formula for the
spectral function $e(\lambda,x,x)$ we need the following
lemma (compare with Lemma~\ref{Singularity of the wave group at time zero lemma}).

\begin{lemma}
\label{Singularity of the wave group at time nonzero lemma pointwise}
Suppose that the point $y\in M$ is nonfocal.
Then for any complex-valued function $\hat\gamma\in C_0^\infty(\mathbb{R})$
with $\operatorname{supp}\hat\gamma\subset(0,+\infty)$ we have
\begin{equation}
\label{Singularity of the wave group at time nonzero lemma pointwise formula}
\mathcal{F}^{-1}_{t\to\lambda}[\hat\gamma(t)\operatorname{tr}u(t,y,y)]=
o(\lambda^{n-1})
\end{equation}
as $\lambda\to+\infty$.
\end{lemma}

\emph{Proof\ }
The result in question is a special case of
Theorem~4.4.9 from \cite{mybook}.~$\square$

\

The following theorem is our main result regarding the spectral function $e(\lambda,x,x)$.

\begin{theorem}
\label{theorem spectral function unmollified two term}
If the point $x\in M$ is nonfocal then
the spectral function $e(\lambda,x,x)$ admits
the two-term asymptotic expansion
(\ref{two-term asymptotic formula for spectral function})
as $\lambda\to+\infty$.
\end{theorem}

\emph{Proof\ }
The result in question is an immediate consequence of
formulae
(\ref{theorem spectral function mollified equation 7}),
Theorem~\ref{theorem spectral function mollified}
and
Lemma~\ref{Singularity of the wave group at time nonzero lemma pointwise}
from the current paper
and Theorem~B.5.1 from \cite{mybook}.~$\square$

\

We now deal with the counting function $N(\lambda)$.

Suppose that we have a Hamiltonian trajectory
$(x^{(j)}(t;y,\eta),\xi^{(j)}(t;y,\eta))$
and a real number $T>0$ such that
$(x^{(j)}(T;y,\eta),\xi^{(j)}(T;y,\eta))=(y,\eta)$.
We will say in this case
that we have a $T$-periodic trajectory originating
from the point $(y,\eta)\in T'M$.

Denote by $(S^*M)^{(j)}$ the unit cosphere bundle,
i.e.~the $(2n-1)$-dimensional surface in the cotangent
bundle defined by the equation $h^{(j)}(y,\eta)=1$.
The unit cosphere bundle is equipped with a natural Lebesgue measure:
the $(2n-1)$-dimensional surface area element on  $(S^*M)^{(j)}$ is
$dy\,d(S^*_yM)^{(j)}$ where $d(S^*_yM)^{(j)}$ is the
$(n-1)$-dimensional surface area
element on the unit cosphere $(S^*_yM)^{(j)}$, see proof of
Lemma~\ref{Singularity of the wave group at time zero lemma}.

Denote by $\Pi^{(j)}$ the set of points in $(S^*M)^{(j)}$
which serve as starting points for periodic trajectories generated by the Hamiltonian
$h^{(j)}$.
It is known, see Lemma 1.3.4 in \cite{mybook}, that the
set $\Pi^{(j)}$ is measurable.

\begin{definition}
\label{definition of nonperiodicity condition 1}
We say that the nonperiodicity condition is fulfilled
if for each
$j=1,\ldots,m^+$ the set $\Pi^{(j)}$ has measure zero.
\end{definition}

We call a $T$-periodic trajectory \emph{absolutely periodic} if
the function
\[
|x^{(j)}(T;y,\eta)-y|^2+|\xi^{(j)}(T;y,\eta)-\eta|^2
\]
has an infinite order zero in the variables $(y,\eta)$, and we denote
by $(\Pi^a)^{(j)}$ the set of points in $(S^*M)^{(j)}$
which serve as starting points for absolutely periodic trajectories
generated by the Hamiltonian $h^{(j)}$.
It is known, see Corollary 1.3.6 in \cite{mybook}, that the
set $(\Pi^a)^{(j)}$ is measurable and,
moreover, the set $\Pi^{(j)}\setminus(\Pi^a)^{(j)}$ has measure zero.
This allows us to reformulate Definition~\ref{definition of nonperiodicity condition 1}
as follows.

\begin{definition}
\label{definition of nonperiodicity condition 2}
We say that the nonperiodicity condition is fulfilled
if for each
$j=1,\ldots,m^+$ the set $(\Pi^a)^{(j)}$ has measure zero.
\end{definition}

In practical applications it is easier to work with
Definition~\ref{definition of nonperiodicity condition 2}
because the set $(\Pi^a)^{(j)}$ is usually much
thinner than the set $\Pi^{(j)}$.

In order to derive a two-term asymptotic formula for the
counting function $N(\lambda)$ we need the following
lemma.

\begin{lemma}
\label{Singularity of the wave group at time nonzero lemma integrated}
Suppose that the nonperiodicity condition is fulfilled.
Then for any complex-valued function $\hat\gamma\in C_0^\infty(\mathbb{R})$
with $\operatorname{supp}\hat\gamma\subset(0,+\infty)$ we have
\begin{equation}
\label{Singularity of the wave group at time nonzero lemma integrated formula}
\int_M
\mathcal{F}^{-1}_{t\to\lambda}[\hat\gamma(t)\operatorname{tr}u(t,y,y)]\,dy=
o(\lambda^{n-1})
\end{equation}
as $\lambda\to+\infty$.
\end{lemma}

\emph{Proof\ }
The result in question is a special case of
Theorem~4.4.1 from \cite{mybook}.~$\square$

\

The following theorem is our main result regarding the counting function $N(\lambda)$.

\begin{theorem}
\label{theorem counting function unmollified two term}
If the nonperiodicity condition is fulfilled then
the counting function $N(\lambda)$ admits
the two-term asymptotic expansion
(\ref{two-term asymptotic formula for counting function})
as $\lambda\to+\infty$.
\end{theorem}

\emph{Proof\ }
The result in question is an immediate consequence of
formulae
(\ref{definition of counting function}),
(\ref{theorem spectral function mollified equation 7}),
Theorem~\ref{theorem spectral function mollified}
and
Lemma~\ref{Singularity of the wave group at time nonzero lemma integrated}
from the current paper
and Theorem~B.5.1 from \cite{mybook}.~$\square$

\section{$\mathrm{U}(m)$ invariance}
\label{U(m) invariance}

We prove in this section that the RHS of formula
(\ref{formula for b(x)})
is invariant under unitary transformations
(\ref{unitary transformation of operator A}),
(\ref{matrix appearing in unitary transformation of operator})
of our operator $A$.
The arguments presented in this section bear some
similarity to those from Section~\ref{U(1) connection},
the main difference being that the unitary matrix-function in question
is now a function on the base manifold~$M$ rather than on $T'M$.

Fix a point $x\in M$ and an index $j$ (index enumerating the eigenvalues
and eigenvectors of the principal symbol) and consider the expression
\begin{multline}
\label{formula for bj(x)}
\int\limits_{h^{(j)}(x,\xi)<1}
\biggl(
[v^{(j)}]^*A_\mathrm{sub}v^{(j)}
\\
-\frac i2
\bigl\{
[v^{(j)}]^*,A_1-h^{(j)},v^{(j)}
\bigr\}
+\frac i{n-1}h^{(j)}\bigl\{[v^{(j)}]^*,v^{(j)}\bigr\}
\biggr)(x,\xi)\,
d\xi\,,
\end{multline}
compare with (\ref{formula for b(x)}).
We will show that this expression
is invariant under the transformation
(\ref{unitary transformation of operator A}),
(\ref{matrix appearing in unitary transformation of operator}).

The transformation
(\ref{unitary transformation of operator A}),
(\ref{matrix appearing in unitary transformation of operator})
induces the following transformation of the principal
and subprincipal symbols of the operator $A$:
\begin{equation}
\label{transformation of the principal symbol}
A_1\mapsto RA_1R^*,
\end{equation}
\begin{equation}
\label{transformation of the subprincipal symbol}
A_\mathrm{sub}\mapsto
RA_\mathrm{sub}R^*
+\frac i2
\left(
R_{x^\alpha}(A_1)_{\xi_\alpha}R^*
-
R(A_1)_{\xi_\alpha}R^*_{x^\alpha}
\right).
\end{equation}
The eigenvalues of the principal symbol remain unchanged,
whereas the eigen\-vectors transform as
\begin{equation}
\label{transformation of the eigenvectors of the principal symbol}
v^{(j)}\mapsto Rv^{(j)}.
\end{equation}
Substituting formulae
(\ref{transformation of the principal symbol})--(\ref{transformation of the eigenvectors of the principal symbol})
into the RHS of
(\ref{formula for bj(x)})
we conclude that the increment of
the expression (\ref{formula for bj(x)}) is
\begin{multline*}
\int\limits_{h^{(j)}(x,\xi)<1}
\biggl(\,
\frac i2[v^{(j)}]^*
\left(
R^*R_{x^\alpha}(A_1)_{\xi_\alpha}-(A_1)_{\xi_\alpha}R^*_{x^\alpha}R
\right)
v^{(j)}
\\
-
\frac i2
\left(
[v^{(j)}]^*R^*_{x^\alpha}R(A_1-h^{(j)})v^{(j)}_{\xi_\alpha}
-
[v^{(j)}_{\xi_\alpha}]^*(A_1-h^{(j)})R^*R_{x^\alpha}v^{(j)}
\right)
\\
+
\frac i{n-1}h^{(j)}
\left(
[v^{(j)}]^*R^*_{x^\alpha}Rv^{(j)}_{\xi_\alpha}
-
[v^{(j)}_{\xi_\alpha}]^*R^*R_{x^\alpha}v^{(j)}
\right)
\biggr)(x,\xi)\,
d\xi\,,
\end{multline*}
which can be rewritten as
\begin{multline*}
-\frac i2\int\limits_{h^{(j)}(x,\xi)<1}
\biggl(
h^{(j)}_{\xi_\alpha}
\left(
[v^{(j)}]^*R^*_{x^\alpha}Rv^{(j)}
-
[v^{(j)}]^*R^*R_{x^\alpha}v^{(j)}
\right)
\\
-\frac 2{n-1}h^{(j)}
\left(
[v^{(j)}]^*R^*_{x^\alpha}Rv^{(j)}_{\xi_\alpha}
-
[v^{(j)}_{\xi_\alpha}]^*R^*R_{x^\alpha}v^{(j)}
\right)
\biggr)(x,\xi)\,
d\xi\,.
\end{multline*}
In view of the identity $R^*R=I$ the above expression can be further simplified,
so that it reads now
\begin{multline}
\label{increment of bj(x) first iteration}
i\int\limits_{h^{(j)}(x,\xi)<1}
\biggl(
h^{(j)}_{\xi_\alpha}[v^{(j)}]^*R^*R_{x^\alpha}v^{(j)}
\\
-\frac1{n-1}h^{(j)}
\left(
[v^{(j)}]^*R^*R_{x^\alpha}v^{(j)}_{\xi_\alpha}
+
[v^{(j)}_{\xi_\alpha}]^*R^*R_{x^\alpha}v^{(j)}
\right)
\biggr)(x,\xi)\,
d\xi\,.
\end{multline}

Denote
$B_\alpha(x):=-iR^*R_{x^\alpha}$
and observe that this set of matrices,
enumerated by the tensor index $\alpha$ running through the values $1,\ldots,n$,
is Hermitian.
Denote also $b_\alpha(x,\xi):=[v^{(j)}]^*B_\alpha v^{(j)}$
and observe that these $b_\alpha$ are positively homogeneous in $\xi$ of degree 0.
Then the expression
(\ref{increment of bj(x) first iteration})
can be rewritten as
\begin{equation*}
\label{increment of bj(x) second iteration}
-
\int\limits_{h^{(j)}(x,\xi)<1}
\left(
h^{(j)}_{\xi_\alpha}\,b_\alpha
-\frac 1{n-1}\,h^{(j)}\,\frac{\partial b_\alpha}{\partial\xi_\alpha}
\right)\!(x,\xi)\,
d\xi\,.
\end{equation*}
Lemma 4.1.4 and formula (1.1.15) from \cite{mybook} tell us
that this expression is zero.

\section{Spectral asymmetry}
\label{Spectral asymmetry}

In this section we deal with the special case when the operator
$A$ is differential (as opposed to pseudodifferential).
Our aim is to examine what happens when we change the sign of the operator.
In other words, we compare the original operator $A$ with the operator
$\tilde A:=-A$. In theoretical physics the transformation
$A\mapsto-A$ would be interpreted as time reversal,
see equation (\ref{dynamic equation most basic}).

It is easy to see that for a differential operator the number $m$
(number of equations in our system) has to be even and that the
principal symbol has to have the same number of positive and negative
eigenvalues.
In the notation of Section~\ref{Main results}
this fact can be expressed as $m=2m^+=2m^-$.

It is also easy to see that the
principal symbols of the two operators, $A$ and $\tilde A$,
and the eigenvalues and eigenvectors of the principal symbols
are related as
\begin{equation}
\label{Spectral asymmetry equation 1}
A_1(x,\xi)=\tilde A_1(x,-\xi),
\end{equation}
\begin{equation}
\label{Spectral asymmetry equation 2}
h^{(j)}(x,\xi)=\tilde h^{(j)}(x,-\xi),
\end{equation}
\begin{equation}
\label{Spectral asymmetry equation 3}
v^{(j)}(x,\xi)=\tilde v^{(j)}(x,-\xi),
\end{equation}
whereas the subprincipal symbols are related as
\begin{equation}
\label{Spectral asymmetry equation 4}
A_\mathrm{sub}(x)=-\tilde A_\mathrm{sub}(x).
\end{equation}

Formulae
(\ref{formula for a(x)}),
(\ref{formula for b(x)}),
(\ref{generalised Poisson bracket on matrix-functions}),
(\ref{Poisson bracket on matrix-functions})
and
(\ref{Spectral asymmetry equation 1})--(\ref{Spectral asymmetry equation 4})
imply
\begin{equation}
\label{Spectral asymmetry equation 5}
a(x)=\tilde a(x),
\qquad
b(x)=-\tilde b(x).
\end{equation}
Substituting (\ref{Spectral asymmetry equation 5}) into
(\ref{a via a(x)}) and (\ref{b via b(x)}) we get
\begin{equation}
\label{Spectral asymmetry equation 6}
a=\tilde a,
\qquad
b=-\tilde b.
\end{equation}

Formulae (\ref{two-term asymptotic formula for counting function})
and (\ref{Spectral asymmetry equation 6}) imply that the spectrum
of a generic first order differential operator is asymmetric about $\lambda=0$.
This phenomenon is known as
\emph{spectral asymmetry}
\cite{atiyah_short_paper,atiyah_part_1,atiyah_part_2,atiyah_part_3}.

If we square our operator $A$ and consider the spectral problem
$A^2v=\lambda^2v$,
then the terms $\pm b\lambda^{n-1}$ cancel
out and the second asymptotic coefficient of the counting function
(as well as the spectral function) of the operator $A^2$ turns to zero.
This is in agreement with the known fact that for an even order semi-bounded
matrix differential operator acting on a manifold without boundary
the second asymptotic coefficient of the counting function is zero, see
Section 6 of \cite{VassilievFuncAn1984} and \cite{SafarovIzv1989}.

\section{Bibliographic review}
\label{Bibliographic review}

To our knowledge, the first publication on the subject
of two-term spectral asymptotics for systems
was Ivrii's 1980 paper \cite{IvriiDoklady1980}
in Section 2 of
which the author stated, without proof, a formula for the second
asymptotic coefficient of the counting function.
In a subsequent 1982 paper \cite{IvriiFuncAn1982}
Ivrii acknowledged that the formula from
\cite{IvriiDoklady1980} was incorrect and gave a new formula, labelled (0.6), followed by a ``proof''.
In his 1984 Springer Lecture Notes \cite{ivrii_springer_lecture_notes}
Ivrii acknowledged on page 226 that both his
previous formulae for the second asymptotic coefficient were
incorrect and stated, without proof, yet another formula.

Roughly at the same time Rozenblyum \cite{grisha} also stated
a formula for the second asymptotic coefficient of the counting function
of a first order system.

The formulae from \cite{IvriiDoklady1980}, \cite{IvriiFuncAn1982} and \cite{grisha}
are fundamentally flawed because they are proportional to the subprincipal
symbol. As our formulae
(\ref{b via b(x)}) and (\ref{formula for b(x)})
show, the second
asymptotic coefficient of the counting function
may be nonzero even when the subprincipal symbol is zero.
This illustrates, yet again, the difference between scalar
operators and systems.

The formula on page 226 of \cite{ivrii_springer_lecture_notes}
gives an algorithm for the calculation of the correction term
designed to take account of the effect
described in the previous paragraph. This algorithm
requires the evaluation of a limit of a complicated expression
involving the integral, over the cotangent bundle,
of the trace of the symbol of the resolvent of the operator $A$
constructed by means of pseudodifferential calculus. This algorithm
was revisited in Ivrii's 1998 book, see formulae (4.3.39) and (4.2.25)
in \cite{ivrii_book}.

The next contributor to the subject was Safarov
who, in his 1989 DSc Thesis~\cite{SafarovDSc}, wrote down a formula
for the second asymptotic coefficient of the counting function
which was ``almost'' correct.
This formula appears in \cite{SafarovDSc} as formula (2.4).
As explained in Section~\ref{Main results},
Safarov lost only the curvature terms
$\,-\frac{ni}{n-1}\int h^{(j)}\{[v^{(j)}]^*,v^{(j)}\}$.
Safarov's DSc Thesis \cite{SafarovDSc} provides arguments which are sufficiently
detailed and we were able to identify the precise point
(page 163) at which the mistake occurred.

In 1998 Nicoll rederived \cite{NicollPhD} Safarov's formula
(\ref{formula for principal symbol of oscillatory integral})
for the principal symbols of the propagator, using a method
slightly different from \cite{SafarovDSc}, but stopped short
of calculating the second asymptotic coefficient
of the counting function.

In 2007 Kamotski and Ruzhansky \cite{kamotski}
performed an analysis of the
propagator of a first order elliptic system based on the
approach of Rozenblyum \cite{grisha}, but stopped short
of calculating the second asymptotic coefficient
of the counting function.

One of the authors of this paper, Vassiliev, considered systems in Section 6 of
his 1984 paper \cite{VassilievFuncAn1984}. However, that paper dealt with systems of a very special type:
differential (as opposed to pseudodifferential) and of even (as opposed to odd) order.
In this case the second asymptotic coefficients
of the counting function and the spectral function vanish, provided the
manifold does not have a boundary.


\begin{thebibliography}{20}

\bibitem{atiyah_short_paper}
M.F.~Atiyah, V.K.~Patodi and I.M.~Singer,
Spectral asymmetry and Riemannian geometry.
{\it Bull.~London Math.~Soc.}~\textbf{5} (1973), 229--234.

\bibitem{atiyah_part_1}
M.F.~Atiyah, V.K.~Patodi and I.M.~Singer,
Spectral asymmetry and Riemannian geometry I.
{\it Math.~Proc.~Camb.~Phil.~Soc.}~\textbf{77} (1975), 43--69.

\bibitem{atiyah_part_2}
M.F.~Atiyah, V.K.~Patodi and I.M.~Singer,
Spectral asymmetry and Riemannian geometry II.
{\it Math.~Proc.~Camb.~Phil.~Soc.}~\textbf{78} (1975), 405--432.

\bibitem{atiyah_part_3}
M.F.~Atiyah, V.K.~Patodi and I.M.~Singer,
Spectral asymmetry and Riemannian geometry III.
{\it Math.~Proc.~Camb.~Phil.~Soc.}~\textbf{79} (1976), 71--99.

\bibitem{DuiGui}
J.J.~Duistermaat and V.W.~Guillemin,
The spectrum of positive elliptic operators and periodic bicharacteristics.
{\it Invent.~Math.}~\textbf{29} (1975), 39--79.

\bibitem{IvriiDoklady1980}
V.~Ivrii,
On the second term of the spectral asymptotics for the
Laplace--Beltrami operator on manifolds with boundary and
for elliptic operators acting in fiberings.
{\it Soviet Mathematics Doklady} \textbf{21} (1980), 300--302.

\bibitem{IvriiFuncAn1982}
V.~Ivrii,
Accurate spectral asymptotics for elliptic operators that act in vector bundles.
{\it Functional Analysis and Its Applications} \textbf{16} (1982), 101--108.

\bibitem{ivrii_springer_lecture_notes}
V.~Ivrii,
{\it Precise spectral asymptotics for elliptic operators
acting in fiberings over manifolds with boundary}.
Lecture Notes in Mathematics \textbf{1100}, Springer, 1984.

\bibitem{ivrii_book}
V.~Ivrii,
{\it Microlocal analysis and precise spectral asymptotics}.
Springer, 1998.

\bibitem{kamotski}
I.~Kamotski and M.~Ruzhansky,
Regularity properties, representation of solutions, and spectral asymptotics
of systems with multiplicities.
{\it Comm.~Partial Differential Equations} \textbf{32} (2007), 1--35.

\bibitem{NicollPhD}
W.J.~Nicoll,
Global oscillatory integrals for solutions of hyperbolic systems.
Ph.D.~thesis, University of Sussex, 1998.

\bibitem{rotman}
R.~Rotman,
The length of a shortest geodesic loop at a point.
{\it J.~Differential Geometry} \textbf{78} (2008), 497--519.

\bibitem{grisha}
G.V.~Rozenblyum,
Spectral asymptotic behavior of elliptic systems.
{\it Journal of Mathematical Sciences} \textbf{21} (1983), 837--850.

\bibitem{sabourau}
S.~Sabourau,
Global and local volume bounds and the shortest geodesic loops.
{\it Communications in Analysis and Geometry} \textbf{12} (2004), 1039--1053.

\bibitem{SafarovIzv1989}
Yu.~Safarov,
Exact asymptotics of the spectrum of a boundary value problem, and periodic billiards.
{\it Mathematics of the USSR - Izvestiya} \textbf{33} (1989), 553--573.

\bibitem{SafarovDSc}
Yu.~Safarov,
Non-classical two-term spectral asymptotics for self-adjoint elliptic operators.
D.Sc.~thesis,
Leningrad Branch of the Steklov Mathematical Institute of the
USSR Academy of Sciences, 1989.
In Russian.

\bibitem{Safarov_Tauberian_Theorems}
Yu.~Safarov,
Fourier Tauberian theorems and applications.
{\it J.~Funct.~Anal.} \textbf{185} (2001), 111--128.

\bibitem{mybook}
Yu.~Safarov and D.~Vassiliev,
{\it The asymptotic distribution of eigenvalues of partial differential operators}.
Amer.~Math.~Soc., Providence (RI), 1997, 1998.

\bibitem{shubin}
M.A.~Shubin,
{\it Pseudodifferential operators and spectral theory}.
Springer, 2001.

\bibitem{VassilievFuncAn1984}
D.~Vassiliev,
Two-term asymptotics of the spectrum of a boundary value problem under an interior reflection of general form.
{\it Functional Analysis and Its Applications} \textbf{18} (1984), 267--277.

\end{thebibliography}
\end{document}